\documentclass[reqno,12pt]{amsart}
\usepackage{graphicx}% Include figure files

\newcommand{\be}{\begin{eqnarray}}
\newcommand{\ee}{\end{eqnarray}}
\newcommand{\bee}{\begin{eqnarray*}}
\newcommand{\eee}{\end{eqnarray*}}

\newcommand{\N}{{\mathbb N}}
\newcommand{\Z}{{\mathbb Z}}

\begin{document}

\large
 
 \title [Carlini and Kepler's equation]
 {Francesco Carlini: Kepler's equation and the asymptotic solution to singular differential equations}
 %{The mathematical contributions by Francesco Carlini: from the asymptotic study of singular ordinary differential equations 
 %to complex roots of the equation $ x ^ x = y $.}
 
 \author {
 Andrea Sacchetti
 }

\address {
Department of Physics, Informatics and Mathematics, University of Modena and Reggio Emilia, Modena, Italy.
}

\email {andrea.sacchetti@unimore.it}

\date {\today}

\thanks {
This work is partially supported by GNFM-INdAM and by the UniMoRe-FIM project ``Modelli e metodi della Fisica Matematica''. \ The author is very grateful to Arrigo Bonisoli, 
Marco Maioli and Annalisa Sacchetti for useful discussions. \ 
The author sincerely thanks the historical archive of the Brera Astronomical Observatory for kindly making available the original material on Francesco Carlini, and   
in particular: copy of the list of Carlini's titles contained in the folders A225/014CAR and A225/015CAR (Historical Archive of the Brera Astronomical Observatory, Fondo Francesco Carlini, cart. 225, fasc. 14 and 15); copy of the letters sent by Schumacher to Carlini coded A.O.B. cart. 137, 1849, 22 and A.O.B. cart. 
137, 1850, 4 (Historical Archive of the Brera Astronomical Observatory, Scientific Correspondence, cart. 137, 1849, letter n. 22 and cart. 137, 1850, letter n. 4); 
copy of Carlini's draft letter contained in the folder A227/002/69CAR (Historical Archive of the Brera Astronomical Observatory, Fondo Francesco Carlini, 
cart. 227, fasc. 2, letter n. 69); and copy of the notes written by Carlini to comment the book of Lacroix contained in the folder A275/001/CAR (Historical Archive of the Brera Astronomical Observatory, Fondo Francesco Carlini, cart. 275, fasc. 1.).
}

\begin {abstract} 
Carlini's career was mainly dedicated to astronomy, but he was also a particularly skilled mathematician. \ 
In this article we collect and analyse his mathematical contributions in detail. \ In particular, in his important Memoir of the year 1817 devoted to Kepler's 
equation he introduced an innovative idea to solve ordinary differential equations with singular perturbations by means of asymptotic expansions. \ In the same Memoir also 
appeared, five years before Laplace's contributions, what is usually called the Laplace limit constant. \ Furthermore, Carlini published other mathematical Memoirs anticipating, 
70 years in advance, the importance of complex branches of the Lambert's special function.

Msc: 01A55, 01A70


\end{abstract}

\maketitle

\section {Introduction} Mathematicians working in the field of asymptotics methods ascribe to Francesco Carlini the merit of having anticipated the WKB techniques by a century. \ The WKB method is currently 
one of the key instruments of the semi-classical techniques widely used for the study of Schr\"odinger's equation, and it is so named after the contributions 
independently given by Gregor Wentzel, Hendrik Anthony Kramers and L\'eon Nicolas Brillouin. \ It consists of connecting the approximate 
solutions across turning points \cite {BM}. \ Actually, approximation formulas, obtained far from turning points, were already independently
used by George Green (1837) and  Joseph Liouville (1837) and they can be also found in an investigation by the Italian astronomer Francesco Carlini (1817). \ For 
instance, in the historical introduction to semiclassical methods given by Nanny Fr\"oman and Per Olof Fr\"oman \cite {Froman}, it is written that

\begin {itemize}
\item [] {``At the beginning of the nineteenth century Carlini treated an important problem in celestial mechanics \cite {Carlini1} \ He considered the motion of a planet 
in an elliptic orbit around the sun, with the perturbations from all other gravitating bodies neglected. \ $\ldots$ \ The problem treated by Carlini was to determine 
the asymptotic behaviour of the coefficients of the sines in this series for large values of the summation index. \ In his treatment of this problem Carlini had 
to investigate a function $s$ of a variable $x$. \ $\ldots $ \ Carlini, who needed a useful approximate formula for this function when its argument is smaller than its 
order $p$, which tends to infinity, showed that $s(x)$ satisfies a linear, second-order differential equation containing the large parameter $p$. \ In this differential equation 
Carlini introduced a new independent variable $y$ by putting $s= \exp \left ( \frac 12 p \int^x y dx \right )$. \ Then he expanded the function $y$ in inverse powers of $p$. \ When 
Carlini introduced this expansion into the differential equation for $y$ and identified terms containing the same power of $1/p$, he obtained recursive formulas which give what is 
now usually called the WKB approximation  $\ldots $ \ Carlini thus automatically achieved in any order of approximation the result that Kramers achieved in the first-order WKB 
approximation by empirically replacing $\ell (\ell +1)$ by $(\ell + 1/2)^2$, where $\ell$ is the orbital angular momentum quantum number. $\ldots$ \ Referring to the historical 
development described here we find it most natural that the so-called WKB approximation (in which $ \ldots $ and Kramers) should be called the Carlini approximation, $\ldots$.''}
\end {itemize}
In fact, the same authors in a later book \cite {FF2} call the WKB (and Jeffreys) techniques as Carlini (JWKB) techniques.

Curiously enough, ample credit for having disseminated Francesco Carlini's Memoir \cite {Carlini1} must be attributed to Carl Gustav Jacob Jacobi, who, after 30 years from the publication of the Memoir in Milan, recovered it, correcting some errors and publishing an amended version in a scientific journal of wide circulation. \ We are not afraid to say that without this intervention by the German mathematician Carlini's fruitful intuitions would have been ignored. \ In this paper we consider in details Carlini's Memoir and Jacobi's contributions and comments; in particular, a large attention is devoted to Carlini's innovation for the solution to the ordinary differential equation (\ref {Car4Bis}) containing a large parameter. 

However, we are not only interested this Memoir.  \ In fact, we have analysed the entire scientific production of the Italian scientist highlighting the most interesting mathematical contributions. \ Indeed, considering all of Carlini's mathematical contributions, it turns out that there are also other topics he dealt with that are worth of mention.

The paper is organised as follows. 

In \S \ref {Sec2} we trace a brief biography of Francesco Carlini and we also discuss his scientific achievements.

In \S \ref {Sec3} we briefly summarise the Memoir by Carlini entitled {\it Ricerche sulla convergenza della serie che serve alla soluzione del problema di Keplero} 
\cite {Carlini1} published in the year 1817, and we comment on Jacobi's contributions and the Carlini-Jacobi correspondence; a more detailed summary of this 
Memoir is given in Appendix \ref {AppA}. \ Since the Memoir concerns celestial mechanics and, in particular, Kepler's equation, we open with a subsection \S 
\ref {Sec3_1} where the mathematical problems connected to Kepler's equation and the series expansion proposed by Lagrange in order to find out the 
solutions to Kepler's equation are discussed. \ We recall also that the same problem was treated by Pierre-Simon Laplace too in a Memoir \cite {Laplace0} published in the year 1823 
and in the {\it Suppl\'ement au $5^\circ$ volume du trait\'e de M\'ecanique C\'eleste}, Paris (1827) \cite {Laplace}. \ In fact, Laplace gave a result concerning  
the convergence of the series for the radius vector without mentioning Carlini's Memoir and using an approach that was judged to be not rigorous by Jacobi.\ Finally, 
we mention the results given by Friedrich Wilhelm Bessel \cite {Bessel}. \ In subsection \S \ref {Sec3_2} we highlight Carlini's Memoir, leaving to the Appendix \ref {AppA} a more 
detailed description. \ In subsection \S \ref {Sec3_3} we comment on the two Memoirs by Jacobi. \ In subsection \S \ref {Sec3_4} the correspondence Carlini-Jacobi is 
considered. \ Finally, in subsection \S \ref {Sec3_5} we consider a Memoir by Carlini again on Kepler's equation entitled {\it Descrizione d'una macchinetta che serve a 
risolvere il problema di Keplero, ossia a trovare l'anomalia eccentrica data l'anomalia media, qualunque sia l'eccentricit\`a} and published in the year 1853, where 
he makes some important comments concerning the errors in his memoir of the year 1817 and Jacobi's contribution; we leave to Appendix \ref {AppB} a more detailed description 
of this Memoir.

In \S \ref {Sec4} we discuss the mathematical achievements by Carlini in the study of the ordinary differential equation (\ref {Car4Bis}) with large parameter $p$ 
comparing his techniques with the methods already developed by the Swiss mathematician Jean Trembley and the French mathematician Sylvestre Fran\c cois Lacroix.

Finally, in \S \ref {Sec5} we make some concluding remarks. \ In particular we emphasize that the constant usually named Laplace limit constant should be renamed as 
Carlini-Laplace limit constant because it has been introduced by Carlini five years before the contributions by Laplace.

In Appendix \ref {AppA} we give a detalied summary of the Memoir entitled {\it Ricerche sulla convergenza della serie che serve alla soluzione del problema di Keplero} with 
a precise description of the errors contained in this Memoir and how they were corrected by Jacobi.

In Appendix \ref {AppB} we briefly comment some other minor mathematical works by Carlini. \ In particular, we discuss also the Memoir by Carlini 
entitled {\it Sopra alcune funzioni esponenziali comprese nella formola $x^{x^n}$} \cite {Esponenziale}, where he gave some contributions to the analysis 
of Lambert's special function.

\section {Francesco Carlini's Biography} \label {Sec2} Let us briefly sketch the biography of Francesco Carlini; it is based on information collected by Giovanni Virginio Schiaparelli
\footnote {G.V. Schiaparelli was born in Savigliano in Piemonte on March 14th, 1835; graduated in Engineering in Turin in 1854, he then studied astronomy at the 
{\it Berlin Observatory} and at the {\it Russian Imperial Observatory of Pulkovo} (St. Petersburg). \ He returned to Italy in 1860 and was named "second astronomer" (typically each Specola had, at that time, two Astronomers  Professors  and 3 students. When Schiaparelli was called to the Brera Observatory the other Astronomer Professor was Francesco Carlini who also had the role of Director of the Observatory) at the {\it Osservatorio Astronomico di Brera} (hereafter, Brera Observatory), in Milan, and, in 1862, Director. \ He died in Milan on July 10th, 1910. \ He is particularly known for his studies of Mars.}
on the occasion of the commemoration session held on December 18th, 1862 at the {\it R. Lombardo Istituto di Scienze e Lettere} \cite {Schiaparelli1} (see also the biography and the references contained in the note by Giulio Cesare Giacobbe \cite {Giacobbe}). \ The events 
concerning Carlini's private life were taken from a manuscript where Carlini himself had recorded them, and from communications provided by his family and by his personal 
doctor, who was a close friend of Carlini for many years. \ Francesco Carlini was born in Milan on January 7th, 1783 by Rosa Minola and Carlo Giuseppe Carlini. \ His 
father was an employee at the Library of Brera since its foundation and he took care of the primary education of his son by introducing him to the first rudiments of Latin; 
then Carlini attended the Brera Gymnasium. \ The astronomical observatory was established in Brera in the year 1760 and it rapidly became, under the guidance of 
Francesco Reggio, Angelo De Cesaris and Barnaba Oriani, one of the most important in Europe \cite {Schiaparelli2}. \ Carlini was immediately attracted to 
the astronomical activities and from the year 1797 he attended the observatory as a volunteer student performing simple calculations on behalf of astronomers. \ In 
1799 he was admitted as fellow at the {\it Specola}, in 1803 he obtained the diploma of Mathematics at the University of Pavia, and in 1804 he was promoted to the rank 
of {\it astronomo soprannumerario}.\footnote {{\it Astronomo soprannumerario} is that student who has particularly distinguished himself and who, when necessary, replaces the Astronomer Professors in the conduct of public lectures or by continuing a series of observations undertaken by them.} \ In 1832, the year in which he married Gabriella Sabatelli, he was appointed Director of the Brera Observatory, succeeding to  
De Cesaris who died the same year, covering this role until 1862, when Schiaparelli became his successor. \ Carlini was of average height and robust
constitution, he lived in good health to his late seventies, as also evidenced by the fact that in 1860 he undertook, on behalf of the government, a 
trip to Spain to observe a total eclipse. \ In the last months of his life symptoms of a severe intestinal disease began to appear; eventually, he died on 
August 29th, 1862.
\begin{center}
\begin{figure}
\includegraphics[height=10cm,width=8cm]{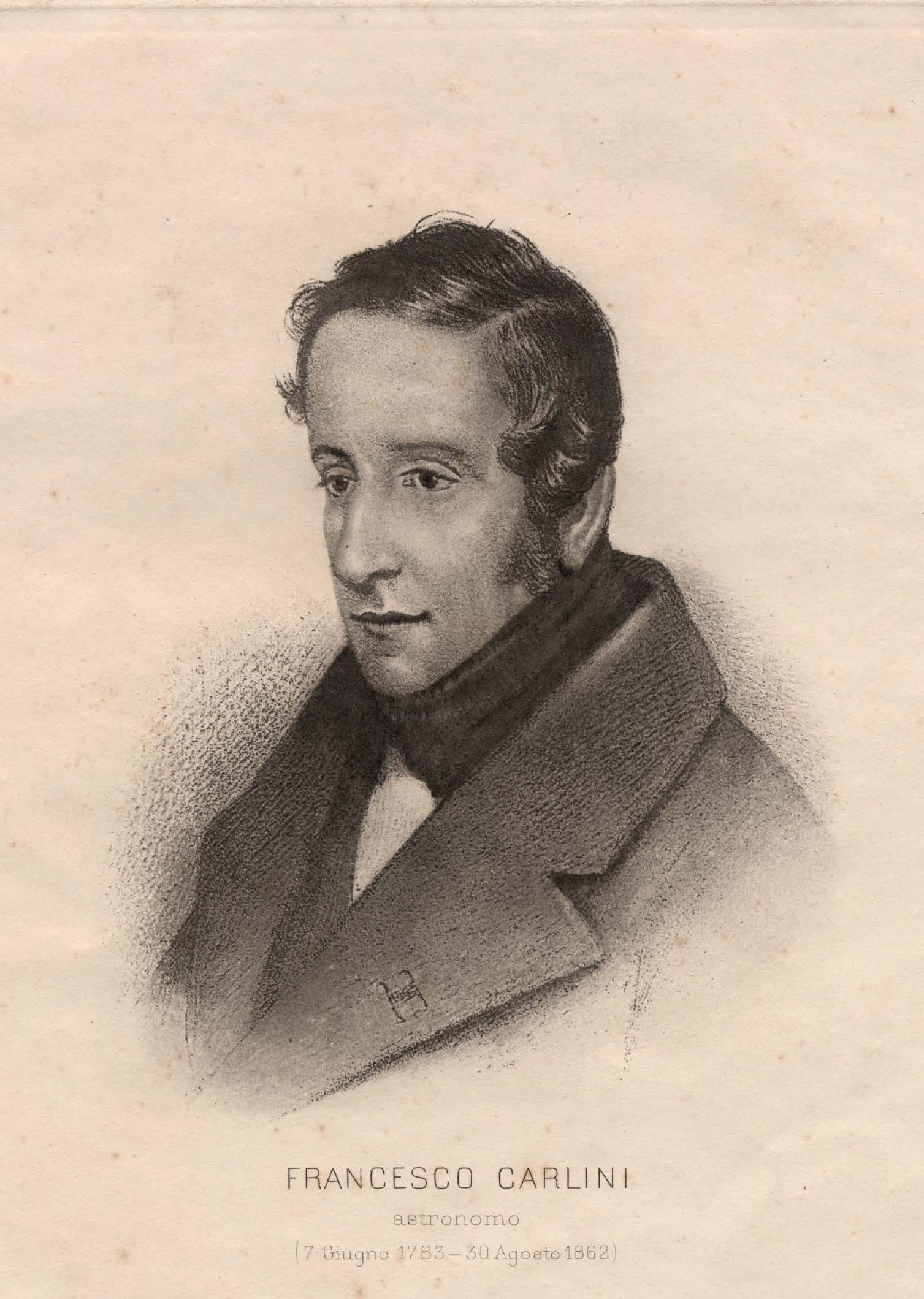}
\caption{\label {Fig1} Image of Francesco Carlini kindly made available by the {\it INAF-Osservatorio di Capodimonte}. \ In fact, the birth and death dates are 
incorrectly reported; indeed Francesco Carlini was born on January 7th, 1783, and he died on June 29th, 1862, according the commemoration Memoir 
\cite {Schiaparelli1}. \ As far as we know, only two pictures of Carlini are available; this one and another one at the Brera Observatory \cite {Froman}.}
\end{figure}
\end{center}

His successes were widely recognised throughout Europe. \ Since 1812 he was a member of the {\it Istituto Nazionale Italiano},  later commuted to
{\it Istituto Lombardo-Veneto } and then to {\it Istituto Lombardo}, and he was its President for a long time. \ He has been a member of over 30 Italian and foreign 
scientific academies.
\begin{center}
\begin{figure}
\includegraphics[height=10cm,width=12cm]{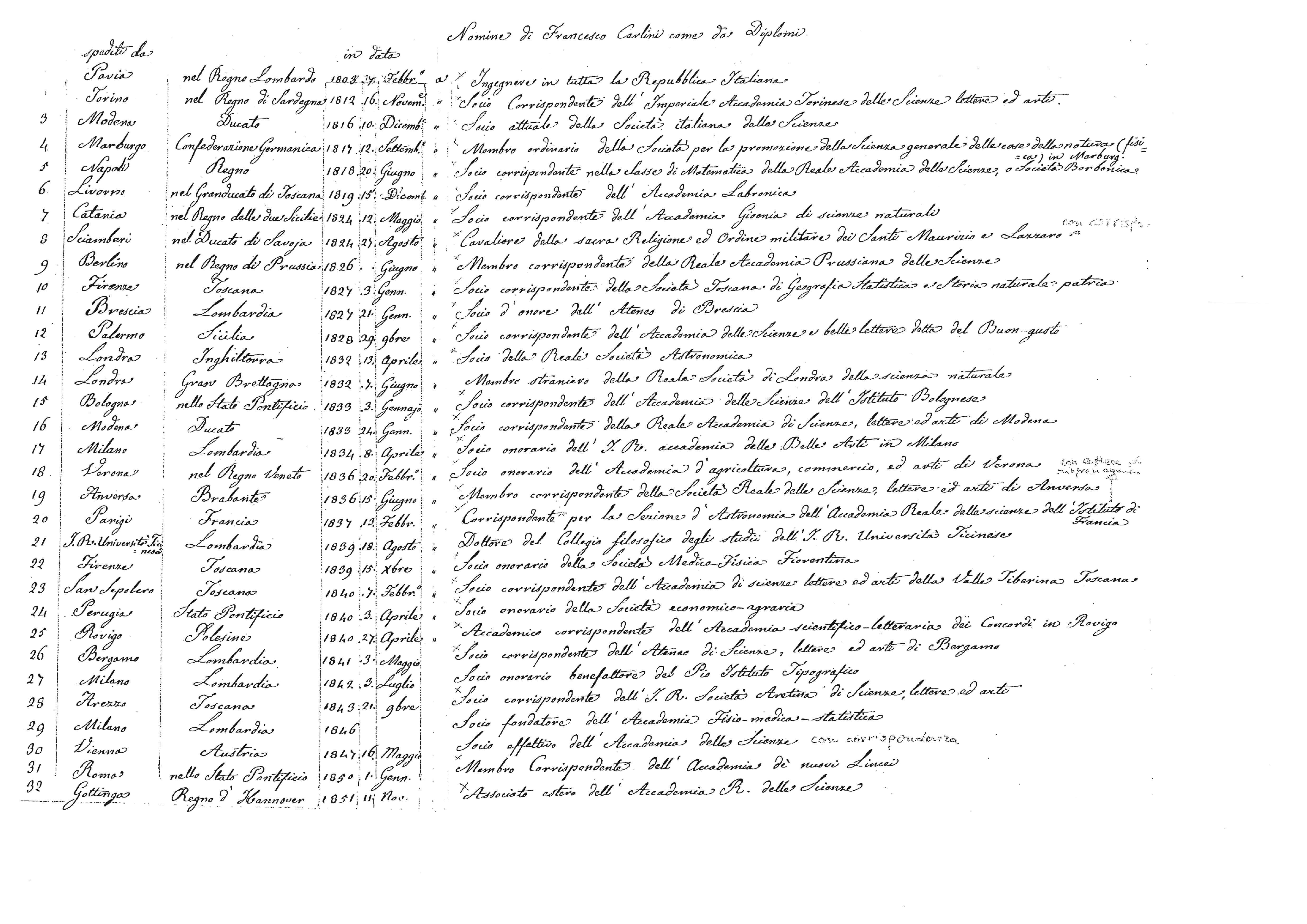}
\caption{\label {Fig2} Image of Francesco Carlini's list of awards made available by the archive of {Brera Observatory}; folders A225/014CAR and A225/015CAR (Historical Archive of the Brera Astronomical Observatory, Fondo Francesco Carlini, cart. 225, fasc. 14 and 15). \ He has been a member of over 30 Italian and foreign 
scientific academies.}
\end{figure}
\end{center}

\subsection {Scientific achievements} On January 1st of the year 1801 the new century opened with the discovery of the asteroid Ceres by Giuseppe Piazzi of the astronomical Observatory of Palermo;
this discovery was soon followed by the contributions of Heinrich Wilhelm Olbers and Karl Ludwig Harding. \ Carlini became interested in these discoveries and in 1803 he 
studied the opposition of the asteroid Pallas, discovered the year before by Olbers. \ Subsequently, from 1814 to 1818, basing his analysis on the works by 
Carl Friedrich Gauss, Carlini produced the tables of the motion of the centre and of the ecliptic reduction of the four main asteroids. \ In this way he began a series of researches 
that occupied the rest of his life, showing his uncommon talent as mathematical researcher, geodesist, astronomer and meteorologist.

During his long career Carlini was in contact with the most important Italian and European astronomers, such as George Biddell Airy and Friedrich Wilhelm Bessel; he 
published at least 145 memoirs listed by Schiaparelli \cite {Schiaparelli1} and so subdivided:

\begin {itemize}

\item [-] Articles published in {\it Effemeridi astronomiche} of Milan (from 1805 to 1863): 67. \ In 1804 Carlini was charged with calculating the astronomical 
ephemerides, which at that time enjoyed great esteem in Europe. 
 
\item [-] Articles published in {\it Memorie dell'Istituto del Regno Lombardo-Veneto} (Milan): 4.

\item [-] Articles published in {\it Memorie dell'Istituto Lombardo} (Milan): 2.

\item [-] Articles published in {\it Atti dell'Istituto Lombardo} (Milan): 1.

\item [-] Articles published in {\it Giornale dell'Istituto Lombardo e Biblioteca Italiana} (Milan): 6.

\item [-] Articles published in  {\it Giornale dell'Istituto Lombardo} (Milan): 3.

\item [-] Articles published in {\it Biblioteca Italiana} (Milan): 19; in 1826 he was appointed co-director, a position he held until 1840.

\item [-] Articles published in {\it Memorie della Societ\`a Italiana delle Scienze} (Modena): 2. 

\item [-] Articles published in {\it Monatliche Correspondenz} by Zach (Gotha): 3.

\item [-] Articles published in {\it Zeitschrift f\"ur Astronomie by Lindenau and Bohnenberger} (T\"ubingen): 1.

\item [-] Articles published in {\it Correspondence Astronomique} by Zach (Genoa): 8.

\item [-] Articles published in {\it Schumacher Astronomische Nachrichten} (Altona): 1; in fact, it is the translation by Jacobi of Carlini's Memoir published 
in 1817 entitled {\it  Ricerche sulla convergenza della serie che serve alla soluzione del problema di Keplero}.

\item [-] Articles published in {\it Gazzetta privilegiata} (Milan): 3.

\item [-] Articles published in {\it Sitzungsberichte der Kaiserlichen Akademie der Wissenschaften 
mathematisch-naturwissenschaftliche Classe}  (Wien): 1.

\item [-] Articles published in {\it Jahrb\"ucher des K.K. Polytechnischen Instituts} (Wien): 1.

\item [-] Articles published in {\it Atti dell'Accademia Fisio-medico-statistica} (Milan): 5.

\item [-] Articles published in {\it Giornale di Fisica e Chimica} by Brugnatelli (Pavia): 2.

\item [-] Articles published in {\it Annuario geografico} by the Earl Annibale Ranuzzi (Bologna): 1.

\item [-] Articles published in {\it Raccolta Scientifica} by the Abbot Palomba (Roma): 3.

\item [-] Memoirs read at Italian Scientific Congresses in different cities: 9. 

\item [-] Works and memoirs printed separately: 3

\end {itemize}

Among the major contributions by Carlini we mention the following ones.

\subsubsection {\bf Improvement of the tables related to the movement of the Sun} When Carlini began his astronomical career, astronomers used 
Jean-Baptiste Delambre's astronomical tables on the movement of the Sun as a reference.   \ Carlini resumed this work 
and, after discovering some minor errors, he undertook a total recalculation of the tables using his own innovative system; he also accompanied them 
with a Memoir entitled {\it Esposizione di un nuovo metodo di costruire la Tavole astronomiche, applicato alle Tavole del Sole} published in Milan in
1810. \ The new method brought a revolution in the calculation methodology adopted until then, introducing advantages that were readily recognised 
by all astronomers, and Carlini's tables acquired a great popularity. \ Later, in 1852, due to the increase in the number of astronomical observations and 
the improvement of the instruments, a new updated edition of the tables was deemed necessary, also following an indirect impulse apparently coming from the 
mathematician and astronomer Friedrich Wilhelm Bessel. \ To give the second edition of the tables a greater precision Carlini considered the effects of the lunar 
motion on the movement of the Earth; in particular, he gave an estimate of the lunar mass equal to $1/80.08$ of the Earth's mass.

\subsubsection {\bf Study of the lunar motion} The study of the lunar motion is one of the classical problems in Celestial Mechanics. \ On this important
problem have been given relevant contributions by many scientists, such as Claudio Tolomeo, Ab\= u al-Waf\= a' and  Tycho Brahe. \ Isaac Newton paved the way for the 
possibility of explaining analytically, through the law of gravitation, most of the characteristics of the lunar motion; however, the complete 
prediction of the lunar movements based on analytical techniques was still undone. \ From 1813 Carlini began a collaboration with the astronomer 
Giovanni Antonio Amedeo Plana to elaborate a complete theory of the lunar motion. \ The aim was to analytically explain all the aberrations, small and large, of the 
satellite and the astronomical observations should only be used to find out the values of the physical constants indispensable for determining the problem. \ This 
approach was much deeper and more complete than the procedures used at that time and it was based on the systematic application of series expansion of functions, 
the use of symbolic notation and the search of algebraic and general solutions, not numerical and particular ones. \ In 1820 Pierre-Simon Laplace proposed to the 
{\it Acad\'emie royale des Sciences} of Paris a prize having as subject the computation, by means of a purely theoretical approach,  of lunar tables exact as those that 
until then had been built with the help of semi-empirical approaches and observations. \ When the announcement of this award came out, Carlini and Plana presented the 
results already obtained and the commission composed by the mathematicians Pierre-Simon Laplace and Sim\'eon-Denis Poisson and the astronomer Johann Karl Burckhardt 
awarded them the first prize, ex aequo with another similar work presented by Marie Charles Th\'eodore de Damoiseau, considerably increasing Carlini's prestige. \ In 
the same year Carlini and Plana jointly published two works, one of which aimed at defending their theory against some objections raised by Laplace. \ 
Subsequently, they applied themselves to edit a work that included their studies on the subject, but their collaboration ended in 1827 due to disagreements on how 
to set up the studies and the entire work was published in three volumes by Plana alone in 1832 (the cooperation between Carlini and Plana, and the detailed description 
of the disagreements, is well explained in the paper by Tagliaferri and Tucci \cite {TT}). \ Carlini himself, using the techniques developed by him for the 
elaboration of the tables concerning the movement of the Sun, composed some tables concerning the motion of the Moon which, although unpublished 
because not considered by him completely satisfactory, were used for several years for the calculation of ephemerides.

\subsubsection {\bf Mathematical Contributions} Carlini was not only an internationally renowned astronomer, but he was also a fine researcher in mathematical analysis as some 
of his contributions show. \ This fact should not surprise us: Celestial Mechanics has always been a rich source of inspiration for mathematical research, and vice 
versa; and many scientists of that time gave important contributions to both disciplines, for instance George Biddell Airy and Friedrich Wilhelm Bessel. \ Carlini 
published several mathematical memoirs such as the one with the title {\it Sopra alcune funzioni esponenziali comprese nella formola $x^{x^n}$}\footnote {``On some 
exponential functions included in the formula $x^{x^n}$.''} \cite {Esponenziale}, published in {\it Memorie dell'I.R. Istituto del Regno Lombardo-Veneto} (Milan, 1819), or 
like that with the title {\it Sulle propriet\`a delle funzioni algebriche conjugate}\footnote {``On the properties of conjugated algebraic functions.''} \cite {FunzioniConiugate} 
presented in 1854 in the {\it Sitzungsberichte der Kaiserlichen Akademie der Wissenschaften 
mathematisch-naturwissenschaftliche Classe} of Wien. \ Carlini did not even disdain to write 
popular mathematical memoirs, such as the brief one entitled {\it Legge dell'inserzione delle foglie nelle piante}\footnote {``Law of leaves insertion in plants.''} \cite {InserzioneFoglie} appeared in {\it Biblioteca Italiana} (Milan, 1837). \ Without a doubt, Carlini's greatest contribution to mathematical analysis 
concerns the study of the elliptical movement of planets contained in the Memoir entitled {\it Ricerche sopra la convergenza della serie che serve a risolvere il problema di 
Keplero}\footnote {``Researches on the convergence of the series that serves to solve the problem of Kepler.''}; in this thorough work Carlini, due to a miscalculation, was led to 
erroneously conclude that the series used to calculate the eccentric anomaly diverges when the eccentricity exceeds a certain threshold. \ This work, published in 1817 in Italian, 
attracted the attention of the scientific community and, in particular, of the mathematician Carl Gustav Jacob Jacobi who published in 1850 a German version correcting the 
mistake. \ This Memoir had a further remarkable merit: for the first time an asymptotic solution to an ordinary differential equation with a \emph {singular} perturbation 
was given. \ We close by recalling a brief Memoir concerning a practical tools useful to numerically solve Kepler's equation. \ It has been published shortly after the Jacobi's 
paper, and probably it has been motivated by it. \ In particular, in this memoir, entitled {\it Descrizione d'una macchinetta che serve a risolvere il problema di Keplero, ossia a 
trovare l'anomalia eccentrica data l'anomalia media, qualunque sia l'eccentricit\`a (con una tavola)}\footnote {``Description of a (small) machine used to solve Kepler's problem, i.e. to 
find the eccentric anomaly given the average anomaly, whatever the eccentricity (with a table)''} and published in the journal of the {\it I.R. 
Istituto Lombardo di Scienze, Lettere ed Arti} (1853) \cite {macchinetta}, the author acknowledges the fact that his paper of the year 1817 contains a mistake and makes some 
interesting statements concerning the rate of convergence of the series which gives the solution to Kepler's equation.

\subsubsection {\bf Practical and geodetic astronomy} Carlini was not only interested in theoretical astronomy, but he was also widely active in practical astronomy and
geodesy. \ He showed how, with the use of a simple telescope equipped with a level bubble and a micrometer, it is possible to measure the time and determine
the latitude; 
%Vedi Schiaparelli, pg.284 in fondo a dx
he applied this method in Spain in 1860 where, by order of the government, he went to observe the total eclipse of July 18th. \ From the year 1821
he participated an international scientific collaboration among the Swiss Confederation, the Kingdom of France, the Kingdom of Sardinia and the 
Lombardo-Veneto Kingdom. \ The scope was to determine the average meridian width from the Atlantic to the Adriatic by taking precise geodetic measurements between Turin and the 
French border. \ These geodetic measurements were then followed by other ones in different places in northern Italy: from the verification of the "Turin degree" to 
the determination of the astronomical position of the cities of Parma and Pavia.

\subsubsection {\bf Meteorology} When in 1832 he succeeded to Angelo De Cesaris in the direction of the Brera Astronomical Observatory, Carlini expanded his studies also to 
meteorology by continuing and improving those researches that had made the Brera Observatory one of the most qualified institutes in this field. \ In 1835 he replaced 
the two daily observations with a more complete and articulated system of three-hourly measurements using more precise instruments. \ He dedicated several memoirs 
to meteorology and collaborated to the foundation of the {\it Societ\`a Meteorologica Lombarda}.

\section {Memoir by Carlini entitled {\it Ricerche sulla convergenza della serie che serve alla soluzione del problema di Keplero}} \label {Sec3}

\subsection {On Kepler's equation} \label {Sec3_1}
If we denote by $\theta$ the eccentric anomaly and by $c$ the eccentricity of the orbit then Kepler's equation 
\be
u = \theta - c \sin \theta \label {Kep0}
\ee
gives the value of the  average anomaly $ u $ by means of a straightforward calculation. \ The value of the radius vector, as function of the eccentric anomaly, is given by 
\be
r = a (1-c \cos \theta ) \label {Kep0radius}
\ee
where $a$ is the length of the semi-major axis. 

However, solving (\ref {Kep0}) for $\theta$ when $u$ and $c$ are given is a much more complicate problem (for a more detailed review we refer to the monographs by 
Colwell \cite {Colwell} and by Bottazzini and Gray \cite {Bottazzini}). \ In fact, equation (\ref {Kep0}) is transcendent in $\theta$ and the solution for this quantity cannot 
be expressed in a finite number of terms \cite {Moulton}. \ Kepler himself believed that his equation (\ref {Kep0}) cannot be solved exactly and many scientists, including Kepler 
himself, proposed methods of obtaining an approximate solution. 

In 1771 Joseph-Louis Lagrange submitted
a note \cite {Lagrange} to the {\it Berlin Academy} regarding the search for the solution to Kepler's equation, by making use of the Lagrange inversion theorem; 
he also came back later to this problem in his book {\it M\'ecanique analytique} \cite {Lag}. \ Lagrange proved that the solution to Kepler's problem can be expressed through a Lagrange series of the form
\be
\theta = u + \sum_{n=1}^\infty A_n \sin (n u)\, , \label {Kep1}
\ee
whose convergence holds true for small values ​​of the eccentricity $ c $. \ Indeed, if $\theta$ is the solution to equation (\ref {Kep0}) then the Lagrange 
Theorem (see \S 7.32 \cite {WW}) states that for any regular function $F(\theta )$ the following power expansion holds true
\bee
F(\theta ) = F(u) + \sum_{n=1}^\infty \frac {c^n}{n!} \frac {d^{n-1}}{du^{n-1}} \left [ F'(u) \sin^n (u) \right ] \, . 
\eee
In particular, by choosing $F(\theta )=\theta $ Lagrange obtained equation (\ref {Kep1}). \ The first terms $A_n$ of the series are quite easy to write; eventually, Lagrange 
gave the explicit expressions of the coefficients $A_n$ for $n=1,2,3$. \ In fact, higher terms are much more complicated than the first ones, and the computational difficulty 
increases very rapidly. \ In addition, a similar expansion also applies to the radius vector:
\be
r = a + \frac 12 a c^2 + a\sum_{n=1}^\infty B_n \cos (n u) \, , \label {Kep1radius}
\ee
for some coefficients $B_n$.

The question also interested Laplace in a Memoir published in 1780. 

In 1817 Carlini resumed Lagrange's research with the aim to calculate the asymptotic behaviour of the coefficients of the power series for large values of the 
index. \ More precisely, he applied Lagrange's Theorem to the function 
\bee
v = F(\theta ) = 2 \mbox {arctan} \left ( \sqrt {\frac {1+c}{1-c}} \tan \left ( \frac {\theta }{2} \right ) \right )
\eee
where $v$ is the true anomaly. \ In such a way he obtained that $v$ and the radius vector can be expanded in power series
\be
v = u + \sum_{p=1}^{\infty} P_p (c) \sin (pu) \label {CarliniSerie1}
\ee
and 
\be 
r = \sum_{p=0}^\infty Q_p (c) \cos (p u) \, . \label {CarliniSerie2}
\ee
He was able to obtain an estimate of the convergence radius of the series (\ref {CarliniSerie1}) and (\ref {CarliniSerie2}). \ Indeed, in his 
ponderous work, which consists of 47 sections spread over 48 pages, Carlini proved the convergence of the Lagrange series (\ref {CarliniSerie1}) for the determination of the true anomaly $ v $ when the eccentricity $ c $ does 
not exceed a value around $ 0.66 $, and the convergence of the Lagrange series (\ref {CarliniSerie2}) for the determination of the radius vector when the 
eccentricity $c$ does not exceed a value around $ 0.62$. \ Unfortunately, Carlini's work contained some errors and inaccuracies; in particular a sign error led him to obtain an 
inaccurate result that was taken up by Jacobi in the year 1849 \cite {Jacobi1}. \ In the following year Jacobi then published the German translation of the 
Memoir by Carlini \cite {Jacobi2} correcting the errors, proving so the convergence of the series (\ref {CarliniSerie1}) and (\ref {CarliniSerie2}) for any 
value of the eccentricity $0\le c <1$. \ Without any doubt the two memoirs by Jacobi had the notable effect of reviving and spreading Carlini's researches, 
which otherwise might have been ignored. \ In fact, just a few years after the works of Jacobi several authors \cite {Scheibner} took up the lines of research set by Carlini.

Laplace, probably ignoring Carlini's work, published in the year 1823 a Memoir and in 1827 a supplement to the 5th volume of his {\it Trait\'e de M\'ecanique C\'eleste} 
\cite {Laplace} where he demonstrated the convergence of the Lagrange series for the determination of the vector radius for values of the eccentricity $c$ less than 
the value $ 0.66195$. \ It seems appropriate to me to dwell here in more detail on Laplace's contributions because, in my opinion, these are not always reported 
in an suitable way. \ In his first work of 1823, which was then taken up without substantial changes in the supplement published in 1827, Laplace initially 
considered a power series for the vector radius (where we choose, for argument's sake, $a=1$) of the form
\be
r = 1 + \frac 12 c^2 - c \cos u - \frac 12 c^2 \cos 2u + \cdots \label {Lap1}
\ee
where the general term of the series has the expression 
\bee
&& \frac {c^p}{1\cdot (p-1)2^{p-1}} \left \{ p^{p-2} \cos pu - p (p-2)^{p-2}\cos [(p-2)u] \right. \\ 
&& \ \ \left. + \frac {p(p-1)}{2!}(p-4)^{p-2} \cos [(p-4)u] + \cdots \right \} \, . 
\eee
Laplace demonstrated the convergence of the series (\ref {Lap1}) when the eccentricity $c$ meets the condition 
\be 
\ln c < \ln \left ( 1+ \sqrt {1+c^2}\right ) - \sqrt {1+c^2}\, , \label {LaplaceLimit0}
\ee
i.e. Laplace estimated that $c$ should be less than $0.66195$. \ Later Laplace demonstrated that the true anomaly can be expressed, as a function of the 
eccentric anomaly, by a series of the type $v=u+\sum_{p=1}^\infty a^{(p)} \sin (pu)$, that is a series of the kind (\ref {CarliniSerie1}), converging for each 
value of the $0\le c <1$ eccentricity. \ Finally, he considered a series development for the vector radius of the type $r=\sum_{p=0}^\infty b^{(p)} \cos (pu)$, that 
is a series of the kind (\ref {CarliniSerie2}), and for this, through arguments considered not very rigorous, he concluded its convergence when the eccentricity 
$c$ is such that 
\be
\frac {e^{\sqrt {1-c^2}}c}{1+\sqrt {1-c^2}} < 1 \label {LaplaceLimit1}
\ee
which is for every $c$. \ Indeed, in his work Laplace made use of a trick in order to get the condition (\ref {LaplaceLimit1}): a passage from real to purely 
imaginary values of the variable $c$. \ This kind of trick was later criticised by Jacobi, and also Watson posed some doubts about the validity of the argument 
used by Laplace (see pg. 7, \cite {W}): ``To anyone is acquainted with the modern theory of asymptotic series, the fallacious character of such reasoning will 
be evident''.

Also  Bessel \cite {Bessel} contributed to the solution to this problem obtaining in 1819 a solution of the form (\ref {Kep1}) with coefficients
$ A_n = \frac 2n J_n (n c) $, introducing the function (later denoted Bessel's function)
\bee
J_n (z) = \frac {1}{\pi} \int_0^\pi \cos \left ( n x - z \sin x \right ) dx \, , \ n \in \Z \, .
\eee
Similarly, he obtained that the coefficients $B_n$ of the series (\ref {Kep1radius}) are given by $B_n= - \frac {2c}{n} J_n' (n c )$. \ We refer to \cite {Dutka} for a more detailed 
discussion of this point. \ According with Taff \cite {Taff} the asymptotic behaviour of $P_p (c)$ and $\frac 2p J_p (pc)$ are the same for large $p$.

In fact, a very large number of analytical and graphical solutions to Kepler's equation have been studied by many authors and a bibliography containing references to 123 
papers on Kepler's equation appeared in the {\it Bulletin Astronomique}, January 1900, p.37-42. \ In this list the work by Carlini is included.

\subsection {Highlights of Carlini's Memoir} \label {Sec3_2}

In his Memoir Carlini\footnote {Contrary to Carlini, who denoted the eccentricity by the letter $ e $, the Neper's constant 
by the letter $i$ and the imaginary unit by $ \sqrt {-1} $ , in this article $ c $ denotes the eccentricity, $ e $ the Neper's constant and $ i $ the imaginary unit 
in accordance with the current notation.} considered the following relation between the true anomaly $v$ and the average anomaly $u$
\bee
v &=& F(u) + \sum_{\ell \in \N \, ,\ \ell \mbox { even}}  \frac {c^\ell}{\ell!} \frac {d^{(\ell-1)} }{du^{(\ell -1)}}\sum_{p=0}^\infty C^{(p)} \cos (pu) + \\ 
&& \ \ + \sum_{\ell \in \N \, ,\ \ell \mbox { odd}} \frac {c^\ell}{\ell !} \frac {d^{(\ell -1)} }{du^{(\ell -1)}}\sum_{p=0}^\infty D^{(p)} \sin (pu) \, , 
\eee
where $F$ is the function defined as 
\bee
F(u ) = 2 \mbox {arctan} \left ( \sqrt {\frac {1+c}{1-c}} \tan \left ( \frac {u }{2} \right ) \right ) \, , 
\eee
$c$ is the eccentricity and an expression of the coefficients $ C^{(p)} $ and $ D^{(p)} $, which actually depend on the index $ \ell $ too, is given. 
\ Therefore, the above series reduces to a series of sine functions $\sin (pu)$ whose coefficients are generically denoted in Carlini's Memoir by $P$. \ In order to get an asymptotic expression of 
the coefficients $P$ for large value of the index Carlini had to solve the second order differential equation 
\be
\frac {d^2s}{dx^2} + \frac {2p+1}{x} \frac {ds}{dx} = \frac {p^2}{\sin^2 (n)} s\, , \label {Car4Bis}
\ee
where $p$ is large parameter. \ So, he faced a completely new problem: the search for the approximate solution of a \emph {singular} differential equation 
(it is singular because $p$ is a parameter which takes large values). \ As a first step he reduced the equation (\ref {Car4Bis}) to a first-order equation 
by setting $ s = e^{\frac 12 p \int y dx} $ 
\be
\frac {dy}{dx} + \frac 12 p y^2 + (2p+1) \frac {y}{x} = \frac {2p}{\sin^2 n} \, , \label {Car4}
\ee
and then he tried to find an approximate value of the solution $ s $ when $ p $ takes very large values. \ His program was clearly stated:
\begin {itemize}
 \item [] {''Wanting to solve by series the value of $ s $ when $ p $ is a very large number, instead of dealing with the integral formula, it is useful to use the 
 differential equation of the first order given in the previous paragraph. \ To get the first term of the development we will start by neglecting the quantities that are 
 not of the order of $ p $. \ To this task it is necessary to examine the order of magnitude of the variables that enter in the calculation.``}.\footnote {\S 35, pg. 35 
 \cite {Carlini1}  
%``Volendo  risolvere in serie il valore di $s$ quando $p$ \`e un numero grandissimo, invece di trattare la formola %integrale, giova servirsi dell'equazione 
%differenziale di primo ordine data nel paragrafo precedente. \ Per avere il primo termine dello svolgimento cominceremo dal trascurare le quantit\`a che non sono 
%dell'ordine di $p$. \ A tale oggetto \`e necessario esaminare l'ordine di grandezza delle variabili che entrano nel calcolo''
}
\end {itemize}
Eventually, by setting
\be
y = Y + \frac {Y'}{p} +  \frac {Y''}{p^2} +  ecc. \label {Eq10BisBis}
\ee
Carlini found a recursive procedure for determining the terms $ Y $, $ Y'$, $ Y'', \ldots $. \ By means of this result Carlini was finally able to prove that 
the above sine series converges when the eccentricity satisfies the condition 
\be 
\ln c < \ln \left ( 1+ \sqrt {1+c^2}\right ) - \sqrt {1+c^2}\, , \label {LaplaceLimit}
\ee
which is satisfied when about $ c < 0.66$ (according with Carlini's calculations). \ In fact, this last result is not rigorous because of some miscalculations, but the procedure 
used is basically correct. \ We must underline, as also Jacobi pointed out, that condition (\ref {LaplaceLimit}) is the same (\ref {LaplaceLimit0}) that 
Laplace \cite {Laplace0,Laplace} obtained some years later for the convergence of the Lagrange series for the determination of the vector radius.

Carlini's procedure is of historical interest because for the first time, as far as we know, a \emph {singular} differential equation is solved by means of a formal series in 
descending powers of the parameter $p$.

\subsection {Contributions by Jacobi} \label {Sec3_3}

Jacobi, thanks to the German Astronomer Johann Franz Encke who introduced him to Carlini's research, was interested in Carlini's Memoir and he published 
two memoirs on this subject in 1849 and 1850. 

The first Memoir \cite {Jacobi1}, entitled {\it \"Uber die ann\"ahernde Bestimmung sehr entfernter Glieder in der Entwickelung der elliptischen Coordinaten nebst 
einer Ausdehnung der Laplaceschen Methode zur Bestimmung der Functionen gerader Zahlen},\footnote {``On the approximate determination of very remote terms in the expansion of 
the elliptical coordinates, in addition to an extension of the Laplace method for determining the functions of large numbers.''} contains the correct solution to the problem of the 
determination of the coefficients denoted by Carlini as $ P '' $ and $ Q $. \ In fact, Carlini gave the expressions of the dominant trends of $P''$ and $Q$ (see, respectively, formulas  (\ref {Eq10Bis}) and (\ref {RefQ}) contained in Appendix A)  but, unfortunately, his result was affected by a mistake. 

In the second Memoir \cite {Jacobi2}, entitled {\it Untersuchungen \"uber die Convergenz der Reihe, durch welche das Kepleresche Problem gel\"ost 
wird. \ Von Franz Carlini},\footnote {``Analysis on the convergence of the series through which Kepler's problem is solved. \ By 
Franz Carlini''} Jacobi published the version of Carlini's Memoir, translating it into German and correcting the errors. \ It must be said that the first 
Memoir of Jacobi, largely based on the Carlini's Memoir, is difficult to read without having the possibility to access Carlini's Memoir, published 30 years 
earlier in Italian which was hard to find. \ Therefore, Jacobi soon realised the need to give a greater diffusion of Carlini's work and he eventually decided to publish it again, this 
time in German and in a widespread journal, correcting the mistakes.

\subsubsection {\bf On the Memoir of Jacobi \cite {Jacobi1} entitled {\it \"Uber die ann\"ahernde Be\-stimmung sehr entfernter Glieder in der Entwickelung der elliptischen 
Coordinaten nebst einer Ausdehnung der Laplaceschen Methode zur Bestimmung der Functionen gerader Zahlen}}

Jacobi acknowledged the great merit of Carlini's contributions as it can be seen from the introduction:

\begin {itemize}

\item [] ''A long time ago, Carlini dealt in a treatise 
{\it Ricerche sulla convergenza della serie che serve alla soluzione del Problema di Keplero}, Milan 1817, $8^\circ$ (48 pages) with the task of determining approximately the 
coefficient of the sine of a very large multiple of the mean anomaly in the expansion of the equation of the centre. \ The means necessary for this purpose also 
allowed him to solve the similar, much easier task with regard to the expansion of the radius vector. \ As far as I know, no other mathematician or astronomer has 
devoted attention to the former task, which is one of the most difficult of its kind. \ In the supplement, compiled from the papers left by Laplace, to the fifth 
volume of {\it M\'ecanique C\'eleste} (1827) Laplace considers the expansions of the {\it true anomaly} and the radius vector when arranging the very same according to the powers 
of eccentricity, and examines the values of the very remote terms of the expansions arranged as such. \ He thereupon considers the now accepted type of arrangement 
of these series according to the cosine or the sine of the multiples of the mean anomaly and also determines for this arrangement the remote terms of the expansion 
of the {\it radius vector}. \ Yet he did not provide a similar approximate determination of the remote terms for the expansion of the {\it true anomaly} arranged 
according to the sine of the multiples of the mean anomaly. \ The mentioned treatise by Carlini could have raised greater interest but, unfortunately, its results are obviously 
incorrect, which is why, so it seems, one may not have believed in pursuing it further, and the same may also have been disregarded by Laplace as well. \ By looking into this 
further, it could have proved to be an instructive work, which is blemished only by some errors irrelevant to the method.``
\end {itemize}
Subsequently Jacobi highlighted the critical point of Carlini's work:
\begin {itemize}

\item [] ''Carlini finds on page 44 (\S 44) of his treatise that the expansion of 
the true anomaly according to the multiples of the mean anomaly stops converging when the eccentricity $e$ is a root of the equation
\bee
\frac {e i ^ {\sqrt {1 + ee}}} {1+ \sqrt {1 + ee}} = 1
\eee
or
\bee
e = 0.66
\eee
where $i$ represents the basis of the natural logarithms (This equation Carlini arrives at by mistake is coincidentally the same as the equation 
whose root may not exceed $e$ according to the later analysis by Laplace ({\it M\'ecanique C\'eleste}, T. V. Suppl. p. 11), when the expansions are supposed to converge 
according to the powers of the eccentricity. \ He further finds at the end of the treatise (\S 47) that the expansion of the radius vector according 
to the multiples of the mean anomaly stops converging when
\bee
e = 0.62.
\eee
This is obviously incorrect. \ As a matter of fact, it is now known, based on strict proof from Dirichlet, Bessel and others, 
%questa è la traduzione di Jacobi, non ci sono referenze
 that all functions within provided limits can 
be expanded into {\it converging} series, which advance according to the cosine or sine of the multiples of their argument when they assume no infinite values 
within the provided limits.``
%{''Carlini findet in S.44 (\S 44) seiner Abhandlung, dass die Entwickelung der wahren Anomalie nach den Vielfachen der mittleren Anomalie zu convergiren 
%aufh\"ort, wenn die Excentricit\"at $e$ eine Wurzel der Gleichung
%
%\bee
%\frac {e i^{\sqrt {1+ee}}}{1+\sqrt {1+ee}} =1
%\eee
%f1
%
%oder 
%f2
%\bee
%e=0.66
%\eee
%
%ist, wo $i$ die Basis der nat\"urlichen Logarithmen bedeutet (Diese Gleichung, auf welche Carlini irrth\"umlicher Weise kommt, ist zuf\"allig dieselbe, wie die 
%Gleichung, deren Wurzel nach der sp\"ateren Untersuchung von Laplace - {\it M\'ecanique C\'eleste}, T.V. Suppl. p.11 - e nicht \"ubersteigen darf, wenn die 
%Entwickelungen nach den Potenzen der Excentricit\"at convergiren sollen). \ Er findet ferner am Schluss der Abhandlung (\S 47), dass die Entwickelung des Radiusvectors 
%nach den Vielfachen der mittleren Anomalie zu convergiren aufh\"ort, wenn
%
%f3
%\bee
%e=0.62 .
%\eee
%
%Dies ist offenbar falsch. \ Denn man weiss jetzt durch strenge Beweise von Dirichlet, Bessel und Anderen, dass sich alle Functionen innerhalb gegebener 
%Grenzen in \emph {convergirende} Reihen, welche nach den Cosinus oder Sinus der Vielfachen ihres Arguments fortschreiten, entwickeln lassen, wenn sie innerhalb 
%der gegebenen Grenzen keine unendlichen Werthe annehmen.``} 
%
\end {itemize}
And finally he concluded the introduction by specifying that:

\begin {itemize}

\item [] ''The errors Carlini must have made could either be due to the method he applied or be only 
accidental oversights of the calculation. \ In the latter case, I could hope to obtain a confirmation, following the path he had taken, of the results I found on the same subject-matter, and I had even fewer concerns looking for these errors - although this proved to be a difficult task - due to the peculiarity of the methods I employed making 
their further confirmation desirable.``
%{''Die Fehler, welche Carlini begangen haben musste, konnten entweder der von ihm angewandten Methode inh\"ariren, oder nur zuf\"allige Versehen der 
%Rechnung sein. \ Im letzteren Falle konnte ich hoffen, auf dem von ihm hetretenen Wege eine Best\"atigung der von mir \"uber denselben Gegenstand gefundenen Resultate 
%zu erhalten, und ich trug um so weniger Bedenken, das, wenn auch m\"uhsame, Aufsuchen dieser Fehler zu unternehmen, da die Eigenth\"umlichkeit der 
%von mir gebrauchten Methoden ihre anderweitige Best\"atigung w\"unschenswerth machte.``}

\end {itemize}

Jacobi, in this memoir, critically reviewed Carlini's results by promptly reporting the errors. \ Initially, he noted that in the expression (\ref {RefQ}) in Appendix \ref {AppA} of the term $ Q $ one must choose the sign minus (remember that the term in (\ref {RefQ}) $ A = - \sqrt {2} $ is a negative number). \ Not only, some other improvements are needed and he finally obtained the correct expression
\bee
Q = -\frac {2\sqrt[4]{1-c^2}}{p\sqrt {p} \sqrt {2\pi}} \left ( \frac {c e^{\sqrt{1-c^2}}}{1+\sqrt{1-c^2}} \right )^p \, . 
\eee
Therefore the series for the radius vector certainly converges when the term
%f14
\bee
\frac {c e^{\sqrt{1-c^2}}}{1+\sqrt{1-c^2}}
\eee
is less than $ 1 $, and because this term is always strictly less than $ 1 $ for any $ c <1 $ it follows the convergence of the series for any value of the 
eccentricity. \ That is, the method proposed by Carlini would rigorously give that result (\ref {LaplaceLimit1}) obtained by Laplace in the 1823 with not 
rigorous arguments. \ In fact, Jacobi did not fail to comment on the result obtained by Laplace concerning the convergence of the radius vector series observing that the methods used by 
Laplace do not have a rigorous justification but a ``divinatory'' character:

\begin {itemize}
\item [] ``The preceding value of $Q$, which is obtained according to these corrections, 
is exactly the same one Laplace finds on page 20, op. cit.. \ The path on which Carlini arrives at the strange expression of the series identified with $s$, has 
complete mathematical sharpness. \ In contrast, the method taken by Laplace lacks a strict substantiation and has a merely prophetic character.``
%{''Der vorstehende Wert von $Q$, welcher nach diesen Berichtigungen erhalten wird, ist genau derselbe, welchen Laplace a.a. () S. 20 findet. Der Weg, 
%auf welchem Carlini zu dem merkw\"urdigen Ausdruck der mit $s$ bezeichneten Reihe gelangt, besitzt die vollkommene mathematische Sch\"arfe. Dagegen entbehrt die 
%von Laplace eingeschlagene Methode einer strengen Begr\"undung und hat einen mehr divinatorischen Charakter.``}
\end {itemize}
while Carlini's method, although affected by calculation errors, is judged to be correct.

Then Jacobi tackled the more difficult problem of calculating the coefficients of the Lagrange series by finding the correct expression for the term $P''$. \ The 
critical point highlighted by Jacobi is that Carlini, due to the sign error in the series (\ref {Car4_4}), finally obtained the wrong expression $ g = \sqrt {1 + c ^ 2} $ in (\ref {Err1}) 
instead of the correct expression found by Jacobi $ g = \sqrt {1-c ^ 2} $. \ Consequently, Carlini incorrectly deduced the convergence of the series when $ c < 0.66 $, 
while, according with the Jacobi analysis, the series turns out to be convergent for any value $c<1$ of the eccentricity. 

In conclusion, Jacobi recognised the importance of Carlini's contributions, as it can be also seen in the
conclusion of his Memoir:
 
\begin {itemize} 
\item [] '' ... However, it can be seen that all the 
essential difficulties of this task were already overcome by Carlini in 1817, and that he was prevented from finding the result himself only by an oversight 
in the signs upon expanding the powers of the sine according to the sine or cosine of the multiples of the angle.``
%{'' ...Das vorstehende Resultat ist eine wesentliche Erg\"anzung der oben angef\"uhrten Untersuchung von Laplace \"uber die Convergenz der beiden 
%verschiedenen Entwickelungsarten der elliptischen Coordinaten.  Ich vermuthe, dass diese L\"ucke Laplace selbst abgehalten hat, diese Arbeit 
%zu ver\"offentlichen. Man sieht aber, dass alle wesentlichen Schwierigkeiten dieser Aufgabe schon von Carlini im Jahre 1817 \"uberwunden waren, 
%und dass er nur durch ein Versehen in den Zeichen bei Entwickelung der Potenzen des Sinus nach den Sinus oder Cosinus der Vielfachen des Winkels 
%verhindert wurde, das Resultat selbst zu finden.``}
\end {itemize}

\subsubsection {\bf On the Memoir of Jacobi \cite {Jacobi2} entitled {\it Untersuchungen \"uber die Convergenz der Reihe, durch welche das Kepleresche Problem gel\"ost 
wird. \ Von Franz Carlini}}

In this Memoir published in the year 1850 Jacobi translated Carlini's work of the year 1817 from Italian to German. \ The Memoir of Jacobi consists of 47 sections and one
Appendix. \ The first 38 sections, with the exception of a few comments, are the literal translation of the Memoir of Carlini. \ Instead, from section 39  
Jacobi had to take Carlini's errors into account and he rewrote the last part of the Memoir in the correct form which is widely different from the Carlini's 
treatment. \ In fact, because of the miscalculations in Carlini's work Jacobi did not simply corrected the sign errors but he had to perform a more complicated 
calculation in order to get the dominant term of the coefficient $P$ of the sine series. \ Eventually, he was able to compute the asymptotic behaviour of $P$ in the limit 
of large index $p$ 
\bee
P = \frac 1p \left ( \alpha e^f \right )^p \left (1 + \frac {4}{3\sqrt {2p \pi }f^3} \right ) 
\eee
where
\bee
f = \sqrt {1-c^2} \ \mbox { and } \ \alpha = \frac {c}{1+\sqrt {1-c^2}} \, . 
\eee
Then, according with the Jacobi's result, the convergence of the sine series follows provided that condition 
\be
\frac {c e^{\sqrt {1-c^2}}}{1+\sqrt {1-c^2}} < 1 \label {Eq11}
\ee
holds true.

We want to emphasise the merits that Jacobi ascribes to Carlini, as explicitly listed in a footnote on the first page:

\begin {itemize}

\item [] ''The original treatise by Carlini has the title: {\it Ricerca sulla convergenza della serie che serve alla soluzione del problema di Keplero. \ Memoria di Francesco 
Carlini. Milano 1817}. \ Although this treatise is blemished by numerous errors, and its results are incorrect, it undisputedly belongs to the most important and significant works 
about the determination of the values of the functions of large numbers due to the method applied therein and the boldness of its composition. \ 
Its reprint therefore seemed desirable. The author expressed the wish - if one intends to make this reproduction of one of his works from his youth - to have 
this done in a German version, a wish complied with herein. \ Here, the necessary improvements have been applied, which had already been suggested in No. 655 of the 
{\it Schumacher Astronomische Nachrichten} (see p. 175 of this volume). \ The parts, which have had to be added in a completely new revision, especially starting at \S 39, have 
been distinguished by quotation marks. \ An annex for the case of very small eccentricities has been added as well.``
%{''Die Originalabhandlung Carlini's f\"uhrt den Titel: {\it Ricerche sulla convergenza della serie che serve alla soluzione del problema di Keplero. \ Memoria 
%di Francesco Carlini. Milano 
%1817}. \ Obgleich diese Abhandlung von zahlreichen Fehlern entstellt ist, und ihre Resultate falsch sind, so geh\"ort dieselbe doch unstreitig wegen der darin angewandten 
%Methode und der K\"uhnheit ihrer Composition zu den wichtigsten und bedeutendsten Arbeiten \"uber die Bestimmung der Werthe der Functionen grosser Zahlen. \ Es schein 
%daher ihr Wiederabdruck w\"unschenswerth. \ Der Herr Verfasser sprach den Wunsch aus, dass, wenn man diese Reproduction einer seiner Jugendarbeiten beabsichtige, 
%dieselbe in einer deutschen Uebertragung geschehen m\"oge, welchem Wunsche man hier nachgekommen ist. \ Es sind hierbei die nothwendigen Verbesserungen angebracht 
%worden, welche schon in N. 665 der Schumacher Astronomischen Nachrichten (cfr. p. 175 dieses Bandes) angedeutet worden sind. \ Die Stellen, welche, besonders von S. 39 an, in 
%ganz neuer Bearbeitung haben hinzugef\"ugt werden m\"ussen, sind durch Anf\"uhrungszeichen unterschieden worden. \ Auch hat man einen auf den Fall sehr kleiner Excentricit\"aten 
%bez\"uglichen Anhang hinzugef\"ugt.``}
%
\end {itemize}

\subsection {Correspondence Jacobi-Carlini} \label {Sec3_4}

Two letters, written in French, sent to Carlini by Heinrich Christian Schumacher, director of the {\it Schumacher Astronomische Nachrichten}, following Jacobi's suggestions, are collected from 
Carlini's archive at the Brera Observatory. \ These letters are catalogued A.O.B. Scientific correspondence cart. 137, 1849, n.24 and A.O.B. Scientific 
correspondence cart. 137, 1850, n.4. \ In the first letter of November 4, 1849, Schumacher informed Carlini that a first Memoir, which confirmed his results (with the 
correction of the sign error) through a new analysis, had just been published by Jacobi. \ He also communicated that Jacobi, in order to complete his studies with a 
second Memoir, was needed to reprint Carlini's Memoir in the {\it Schumacher Astronomische Nachrichten} so that Carlini's research would have a deserved diffusion. \ In the second 
letter, sent by Schumacher on January 11th, 1850, Carlini was informed that Jacobi had undertaken the translation of his Memoir.

The folder A227/002/45 of the Brera archive contains a draft of a letter, written in Italian language. \ From the content of the letter we think that it could be the 
answer by Carlini to the first letter. \ 
%received from Jacobi; therefore, we assume that this draft was then translated into French (that was the conventionally used language in communications) and sent to Jacobi. 
The content of the letter is the following one: 
\begin {itemize}
\item []{''Dear Sir, I give you my warmest thanks for the precious gift of your learned pamphlet on Kepler's problem, as well as for the mention that in it you have 
had the goodness to make an old calculation of mine on this same subject. \ I am also pleased to have achieved the purpose I set for myself, that is to 
draw the attention of talented mathematicians of Europe to the usefulness that can arise from considering in complicated form the coefficients given by the development of 
the equation of the centre in series of the sinuses of the multiples of the anomaly. \ The introduction of the integrals actually gives, as Mr. Poisson very well 
observes, the simplest and most natural solution to such a problem. \ I was also astonished to see how this author, by publishing his formula as something new, was 
able to ignore the research published since 1818 by the famous Bessel in the journal entitled Zeitschrift f\"ur Astronomie etc. herausgegeben.``} 
\footnote {''Pregiatissimo Signore,
le faccio i miei pi\`u vivi ringraziamenti tanto pel pregiato dono del suo dotto opuscolo sul problema di Keplero, quanto per la menzione che in esso Ella ha avuto 
la bont\`a di fare d'un mio antico abbozzo di calcolo sopra questo medesimo argomento. \ Io mi compiaccio poi di avere ottenuto lo scopo che allora mi ero principalmente 
prefisso, cio\'e di richiamare l'attenzione di valenti matematici d'Europa sopra l'utilit\`a che pu\`o nascere dal considerare sotto forma complessa i coefficienti 
dati dallo svolgimento dell'equazione del centro in serie dei seni dei multipli dell'anomalia. \ L'introduzione degli integrali definiti d\`a realmente, come assai 
bene osserva il Sig. Poisson, la soluzione la pi\`u semplice e la pi\`u naturale d'un tal problema. \ Mi ha fatto per altro meraviglia il vedere come questo autore, 
pubblicando come cosa nuova la sua formula, abbia potuto ignorare le ricerche pubblicate fin dal 1818 dal celebre Bessel nel giornale intitolato 
Zeitschrift f\"ur Astronomie etc. herausgegeben.''}
\end {itemize}
We cannot ignore that there is also a slight criticism of Poisson for ignoring Bessel's research in his studies.

\begin{center}
\begin{figure}
\includegraphics[height=8cm,width=6cm]{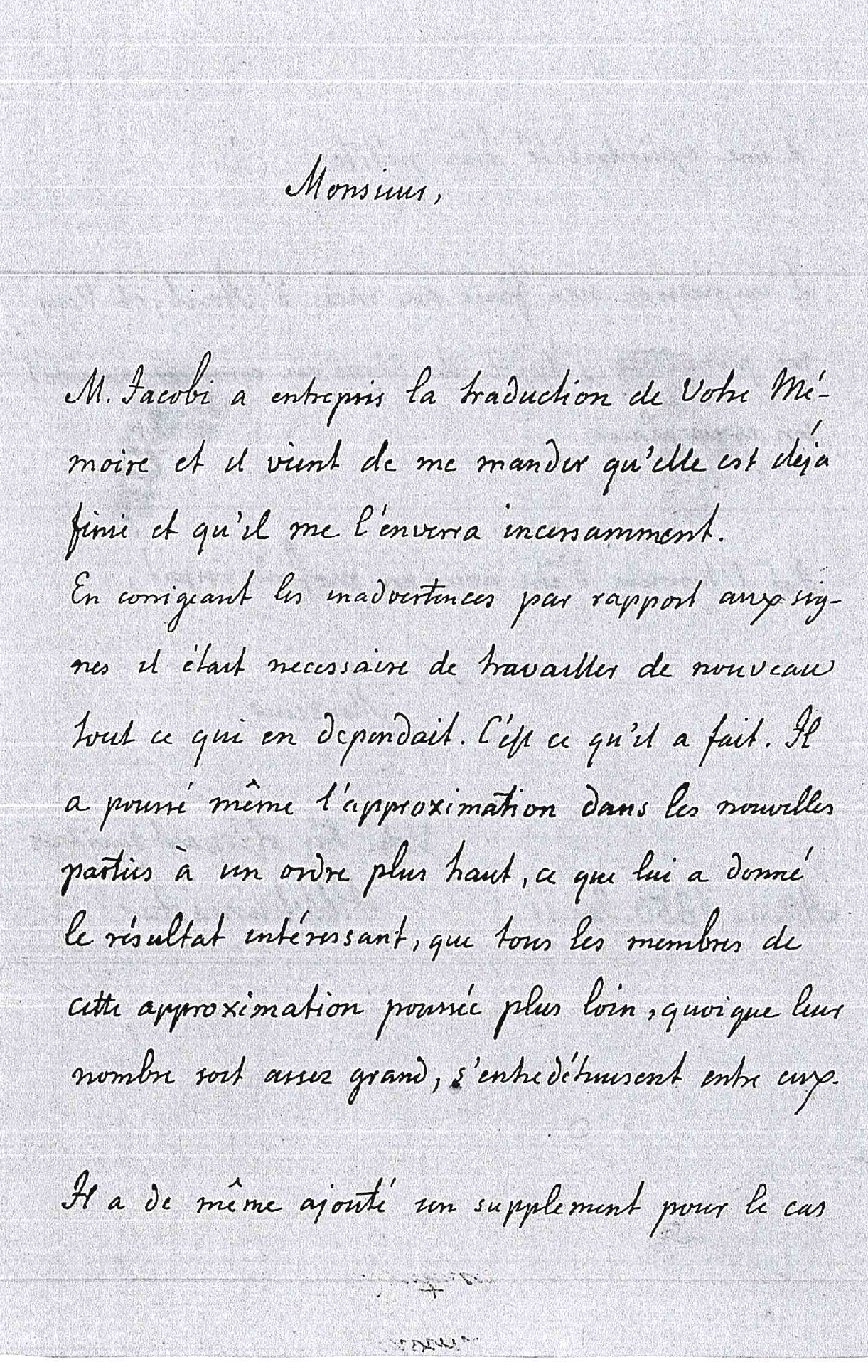}
\includegraphics[height=8cm,width=6cm]{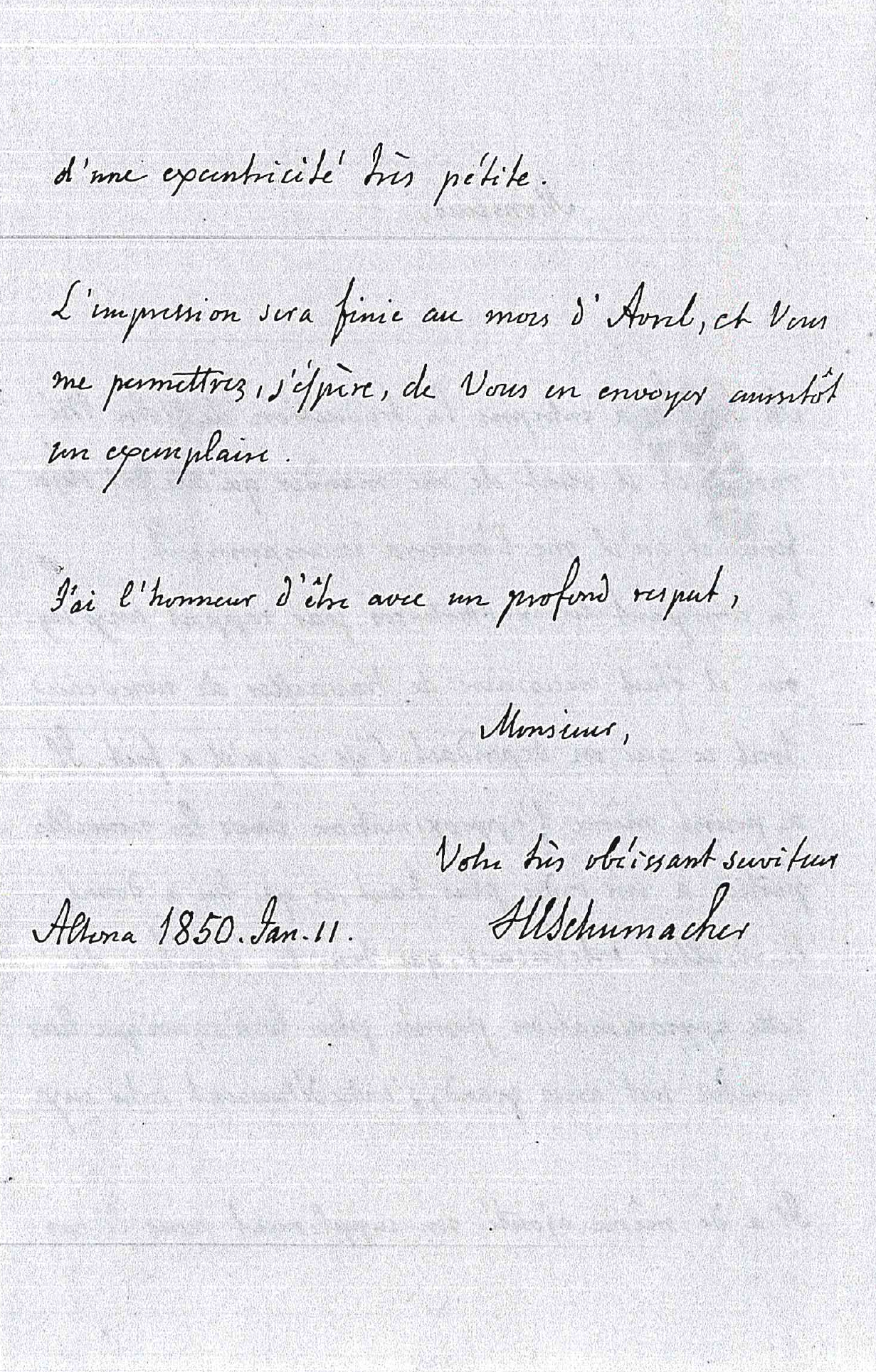}
\caption{\label {Fig4} Copy of the second letter from Schumacker to Carlini; kindly made available by the archive of the {Brera Observatory} (A.O.B. cart. 137, 1850, 4 - Historical Archive of the Brera Astronomical Observatory, Scientific Correspondence, cart. 137, 1850, letter n. 4).}
\end{figure}
\end{center}

\subsection {On a brief Memoir by Carlini concerning a ``macchinetta'' %\footnote {``small machine''} 
that serves to solve Kepler's equation} \label {Sec3_5}
On August 5, 1852, Carlini discussed a report on occasion of the meeting of the Istituto Lombardo Veneto. \ He proposed a small machine useful for the numerical calculation 
of the solution to Kepler's equation (\ref {Kep0}). \ We don't enter here in the details of the memoir, postponing this discussion to the Appendix \ref {Macchi}; however, we cannot 
ignore this important sentence of the author: 

\begin {itemize}
 \item []{``In a paper I published in my youth on this problem... ...I believed I could prove that they cease to converge when eccentricity passes a certain limit. \ But the famous German mathematician Jacobi...exposed, among others, a miscalculation made by me, showing that the series, as given by Lagrange, Oriani and others, is convergent whatever the eccentricity. \ I therefore make my retraction here; but I must warn at the same time that the 
convergence of the series can, beyond certain limits, be so slow that the search, with the help of it, for a single equation of the centre requires many months of work.}\footnote 
{Pg. 440 \cite {macchinetta} ``In uno scritto ch'io pubblicai in mia giovent\`u su tal problema, ... , avendo creduto di poter dimostrare che esse cessano d'essere convergenti quando l'eccentricit\`a 
passa un limite determinato. \ Ma il celebre matematico tedesco Jacobi ... fece avvertire, fra gli altri, uno sbaglio di calcolo da me commesso, corretto il quale si 
dimostra che la serie, come \`e data dal Lagrange, dall'Oriani e da altri, \`e convergente qualunque sia 
l'eccentricit\`a. \ Io faccio dunque qui la mia ritrattazione; ma 
devo avvertire nello stesso tempo che la convergenza della serie pu\`o, oltre certi limiti, essere tanto lenta che la ricerca, col mezzo di essa, d'una sola equazione del 
centro richieda molti mesi di lavoro.''} 
\end {itemize}

In this Memoir, written shortly after Jacobi's publications, the author confirmed the fact that his article published in 1817 contained 
an error that falsified the results and also pointed out that the Lagrange series, which Jacobi proved to be converging for any value of eccentricity, actually converges too 
slowly when the eccentricity is not small enough.

\section {Comment to the "semiclassical" study of the differential equation (\ref {Car4Bis})} \label {Sec4}

\begin{center}
\begin{figure}
\includegraphics[height=12cm,width=10cm]{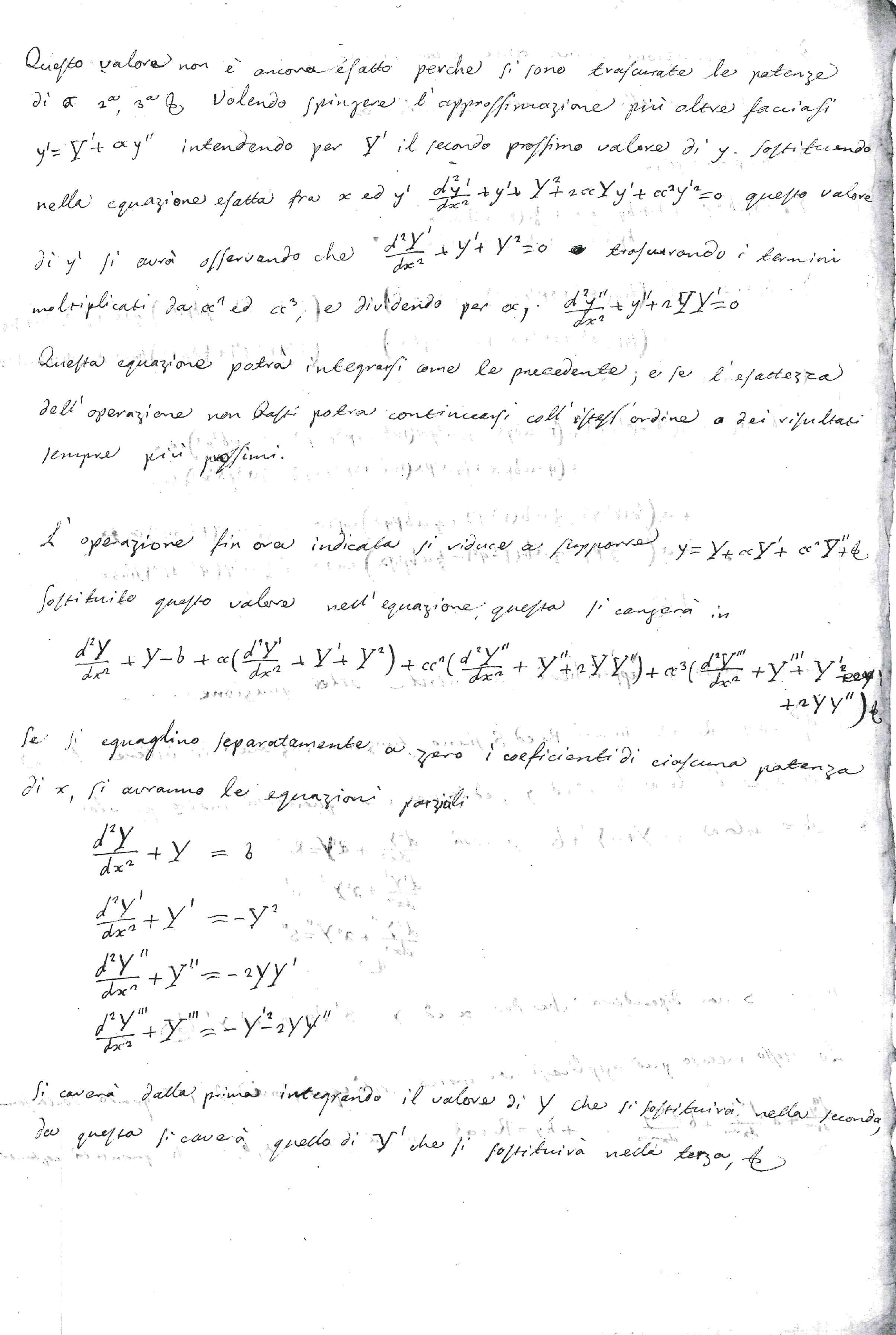}
\caption{\label {Fig3} Carlini's note of the study of the equation (\ref {Eq8Bis}) from the Lacroix's book; kindly made available by the archive of the {Brera Observatory}  
(folder A275/001/CAR - Historical Archive of the Brera Astronomical Observatory, Fondo Francesco Carlini, cart. 275, fasc. 1.). \ One can clearly see the formal 
expansion (\ref {circ}) and the iterative scheme (\ref {circBis}).}
\end{figure}
\end{center}

Now, let's focus our attention to the method used by Carlini to give an approximate expression of the solution to the differential equation (\ref {Car4Bis}). \ Carlini 
ingeniously observed that the solution can be expressed through a (formal) series of functions of the form (\ref {Eq10BisBis}), where $p$ is a parameter which takes 
large values. \ This idea was also independently developed later in the year 1837 by Joseph Liouville and by George Green; and is currently one of the key instruments 
of the semi-classical techniques widely used for the study of Schr\"odinger's equation, under the name of WKB method (from Gregor Wentzel, Hendrik Anthony Kramers and 
L\'eon Nicolas Brillouin). \ In his work Schlissel \cite {Schlissel} analyses the historical path that has characterised the search for asymptotic solutions 
to ordinary linear differential equations; he claims that the very first contributions to this theory were the ones given in the Memoir published by Carlini in the year 
1817. \ In fact, in the work of Schlissel the description of the equation studied by Carlini is not really equation (\ref {Car4Bis}), but a different equation. \ Also 
N. Fr\"oman and P.O. Fr\"oman recognise Carlini's priority in having developed these techniques.

In my opinion this merit must be actually given to Carlini, having however in mind that at that time some methods for the approximate study of differential 
equations were already used. 

In fact, it is certain that Carlini was trained, as a mathematician, studying the works of Sylvestre Fran\c cois Lacroix; this is an ascertainable fact from the 
analysis of Carlini's manuscripts collected in the Carlini's archive at the Brera Observatory. \ Indeed, Lacroix's Treatise played an essential role in the 
transmission of the achievements of the 17th and 18th centuries concerning the series method (see \S 2.3.3 ``Transmission par Lacroix'' \cite {Tournes}). \ In the first 
edition of the year 1798 of the {\it Trait\'e du calcul diff\'erentiel et du calcul int\'egral} \cite {Lac1} Lacroix 
spent two sections of the Chapter III for the explanation of methods to solve by approximations differential equations of the first and second order:
\begin {itemize}

\item [-] {\it M\'ethodes pour r\'esoudre par approximation les \'equations diff\'erentielles du premier ordre};

\item [-] {\it  M\'ethodes pour r\'esoudre par approximation les \'equations diff\'erentielles du second ordre}.

\end {itemize}
As one can easily check it turns out that in these two sections different methods are proposed but no methods similar to those used by Carlini are mentioned, even if 
sometimes the idea to look for a formal solution written by means of power series is applied. 

Conversely, in the second edition \cite {Lac2} published in 1814 Lacroix completely rewrote these two sections in a dedicated chapter (chapter VI) with three expanded 
sections:
\begin {itemize}

\item [-] {\it Des M\'ethodes pour int\'egrer par approximation les \'equations diff\'erentielles};

\item [-] {\it D\'eveloppemens en fractions continues};

\item [-] {\it Usage des \'equations diff\'erentielles du pemier degr\'e pour int\'egrer par approximation}.

\end {itemize}
In particular, in this chapter equations of the form 
\be
\frac {d^2 y}{dx^2} + y + \alpha y^2 = b \label {Eq8Bis}
\ee
are also studied, where $ \alpha $ designates an infinitesimal amount. \ With a modern language we say that the nonlinear term $\alpha y^2$ represents a perturbation and, because $\alpha$ 
takes small values, that (\ref {Eq8Bis}) is an example of \emph {regular} perturbation of the differential equation $\frac {d^2 y}{dx^2} + y = b$ which admits an 
explicit solution. \ In sections 672 and 673 Lacroix, referring to the work by Jean Trembley \cite {Tre} where he studied an equation proposed by Laplace, looked 
for a solution of the form $ y = Y + \alpha Y '$ or, more generally, through the formal development 
\be
y = Y + \alpha Y' + \alpha ^ 2 Y'' + \alpha ^ 3 Y'''+ etc. \label {circ}
\ee
where $ Y $, $ Y '$, $ Y''$ and $ Y''' $, etc., are solutions to the equations
\be
\frac {d^2 Y}{dx^2} + Y  &=& b \nonumber \\
\frac {d^2 Y'}{dx^2} + Y'  &=& -Y^2 \nonumber \\
\frac {d^2 Y''}{dx^2} + Y''  &=& -2Y Y' \label {circBis} \\
\frac {d^2 Y'''}{dx^2} + Y'''  &=& -(Y')^2 -2Y Y'' \nonumber \\
ecc. && \nonumber
\ee
obtained by inserting (\ref {circ}) in (\ref {Eq8Bis}) and equating the terms with the same order of the small parameter $\alpha$. \ Therefore, the idea to formally 
solve a differential equation with a \emph {regular} perturbation by means of a series of functions with ascending powers coefficients $\alpha^{+n}$, $n=0,1,2,3,\ldots $,  
was introduced before the Carlini's Memoir and it was evidently known to Carlini as it appears from Figure \ref {Fig3}.

Carlini's great merit was to realise that this procedure could also be applied to the case of \emph {singular} perturbations of the kind  (\ref {Car4Bis}) where the 
parameter $ p $ is not small, but rather large, and where the formal solution could be written by means of a series (\ref {Eq10BisBis}) of functions with descending 
powers coefficients $ p ^ {- n} $, $n=0,1,2,3,\ldots $.

\section {Conclusion} \label {Sec5}

In his rich and long scientific career Francesco Carlini showed a remarkable talent and interest in the study of Mathematical Sciences. \ No doubt that his main accomplishment in this field was is the one related to the study of the solution of Kepler's equation through series expansion, where the most innovative element is the introduction 
of the formal expansion method in order to search for solutions to ordinary differential equations with \emph {singular} perturbation. \
His great merit was to have realised that the procedure, previously developed by Trembley and Lacroix in the case in which the perturbative parameter is small, 
could also be adapted to the more difficult case in which the perturbative parameter is large. \ This was the starting point from which the asymptotic theory developed up 
to the WKB method.

Another contribution, unfortunately unknown, given by Carlini was the analysis of the function $x^{x}$, with particular regard to the discussion of the complex-valued solutions to the equation 
$x^x=y$ \cite {Esponenziale}. \ The study of the roots in the complex plan became of great interest in the second half of the 1800s and the considerations made by Carlini at the 
beginning of the 1800s were not taken into due consideration. \ With the present article we also want to recall this unacknowledged result to the scientific community.

We want to close with two comments. 

Carlini's article of 1817 was spoiled by an error which falsified its conclusions. \ But the setting and the underlying ideas were 
very interesting and Jacobi had the great merit of recognising its value and giving it a correct version. \ In this way Carlini's work became much more widespread. \ We can 
therefore perhaps conclude that it was a ``lucky'' event that Carlini made that mistake, because otherwise his merit as mathematician would perhaps not have been so widely 
recognised.
 
In mathematical literature ``Laplace limit'' is understood as that value of the eccentricity for which a solution to Kepler's equation, in 
terms of a power series, converges (see Arnold \cite {A} p.85, see also \S 4.8 Laplace Limit Constant 
in the book \cite {Finch}). \ It is approximately $0.66274$ and it is identified with the root of the equation
\be
\frac {c e^{\sqrt {1+c^2}}}{1+\sqrt {1+c^2}}=1 \, . \label {Arnold}
\ee
As far as we know, Carlini was the first author who solved this equation obtaining that the value of the solution is about $0.66$; later Laplace obtained 
for $c$ the value $0.66195$ and finally Cauchy in the year 1829 gave the precise value $0.66274$ \cite {Colwell}. \ In fact, equation (\ref {Arnold}) appears, for the 
first time, in Carlini's Memory published in the year 1817 where he stated that the sine series converges 
for the value of $c$ satisfying condition $\frac {c e^{\sqrt {1+c^2}}}{1+\sqrt {1+c^2}}<1$. \ It is also worth remembering that Laplace \cite {Laplace0,Laplace} discussed 
equation (\ref {Arnold}) concerning the 
convergence of the series used to solve the equation for the radius vector, not Kepler's equation (see \S 100 by \cite {Moulton}). 

\appendix

\section {Summary of the Memoir entitled {\it Ricerche sulla convergenza della serie che serve alla soluzione del problema di Keplero}} \label {AppA}

Here, a summary of the Memoir by Carlini, which consists of 47 sections, follows. 

In the first 4 sections, Carlini describes the reasons that led him to study the problem, which consists in solving Kepler's equation (\ref {Kep0}). \ If one  
denotes by $ v $ the true anomaly then it is linked with the eccentric anomaly $\theta$ by the relation
\bee
v = F(\theta ) = 2 \mbox {arctan} \left ( \sqrt {\frac {1+c}{1-c}} \tan \left ( \frac {\theta }{2} \right ) \right ) 
\eee
and, using the Lagrange Theorem, $ v $ can be eventually written through the series depending from the average anomaly $ u $
\be
v 
%&=& 
%F(u) + \frac {c}{1} \sin (u) F' (u) + \frac {c^2}{1\cdot 2} \frac {d \sin^2 (u) F'(u)}{du} +  \frac {c^3}{1\cdot 2\cdot 3} \frac {d^2 \sin^3 (u) F'(u)}{du^2}
%+ ecc. \nonumber \\
%&=& 
= F(u) + \sum_{\ell =1}^{\infty}  \frac {c^\ell}{\ell !} \frac {d^{(\ell -1)} \sin^\ell (u) F'(u)}{du^{(\ell -1)}} \, , \label {Car2Bis}
\ee
where
\bee
F'(u) = \frac {\sqrt {1-c^2}}{1-c\cos (u)} = \frac {{1-\alpha^2}}{1-2\alpha \cos (u) + \alpha^2} = 1 + \sum_{m=1}^{\infty} 2 \alpha^m \cos (mu ) \, . 
\eee
For argument's sake the units of measure are chosen such that the major semi-axis has unit length, $ f = \sqrt {1-c^2} $ is the length of the minor semi-axis and
$ \alpha = \frac {c} {1 + f} $. 

Starting from section 5 and up to section 16 he considers the terms of the series (\ref {Car2Bis}) with even index $ \ell = 2q $; from the trigonometric relation
\bee
\sin^{2q} (u) = \sum_{j=0}^{q} B^{(2j)} \cos (2ju)
\eee
he finds that the expression $ F '(u) \sin^{2q} (u) = \sum_{p = 0}^\infty C^{(p)} \cos (pu) $ can be written as a trigonometric series of cosine functions. 

In particular, in section 5 he starts considering the case $ p \ge 2q $ finding that 
\bee
C^{(p)} = 2 \alpha^p \sin^{2q} (n) \, , 
\eee
being $ n $ a complex number defined by the relation $\alpha = e^{ in }$. 

When the index $ p $ is less than $ 2q $ then a simple substitution as done above does not work and, in sections 6 and 7, he considers the case $ p = 0 $ finding, 
through a sequences of quite complex but tedious calculations, the following relation
\bee
C^{(0)} = 2 q i B^{(0)} \sin^{2q} (n) \int \frac {dn}{\sin^{2q} (n)} + C
\eee
where the integration constant $ C $ is determined such that  $ C^{(0)} = B^{(0)} $ when $ \alpha = 0$. 

Furthermore, in section 8 he obtains the following asymptotic behaviour when $ q $ is very large: \footnote {Actually, the correct asymptotic behaviour, as indicated by Jacobi, is $C^{(0)} = B^{(0)} \sqrt {1-c^2} \left [ 1 + \frac {c^2}{2q} + ecc. \right ]$. \ In fact, this is a minor mistake that does not affect the validity of the construction.}
\bee
C^{(0)} = B^{(0)} \sqrt {1-c^2} \left [ 1 + \frac {c^2}{q} + ecc. \right ] \, . 
\eee

From section 9 to section 14 he studies the intermediate case $ 0 <p <2q $ considering, at the first instance, the case in which $ p = 2r $ is an even index. \ After a 
series of cumbersome calculations he finally finds that for large $ q $ and $ r $ the following asymptotic behaviour 
holds true
\bee
C^{(2r)} = - \frac {4 i B^{(2r)} (q-r)(q+r+1)}{(2q+1)^2 \mbox {cot} (n) + (2r+1)^2 \tan (n) } + ecc. 
\eee

In section 15 he then considers the case in which $ 0 <p = (2r + 1) <2q $ is an odd index obtaining a similar result. 

In section 16 an iterative relation between the coefficients $ C^{(p)} $, by which it is possible to determine $ C^{(1)} $ starting from $ C^{(0)} $, $ C^{(2)} $ starting 
from $ C^{(1)} $, and so on, is given. 

Starting from section 17 and up to section 19 he considers the terms of the series (\ref {Car2Bis}) with odd index $ \ell = 2q + 1 $ and, by making use of 
the trigonometric relation
\bee
\sin^{2q+1} (u) = \sum_{j=0}^{q} B^{(2j+1)} \sin \left [ (2j+1)u\right ]\, , 
\eee
he finds that the expression of $ F '(u) \sin^{2q + 1} (u) = \sum_{p = 1}^\infty D^{(p)} \sin (pu) $ can be written as a trigonometric series of sine functions 
with coefficients $ D^{(1)} $, $ D^{(2)}, \ldots $. \ Similarly to the calculation of the coefficients $ C^{(p)} $, the coefficients $ D^{(p)} $ are determined
considering at first the case in which $ p \ge 2q + 1 $, subsequently the case $ p = 0 $ and finally the intermediate case $ 0 <p <2q + 1 $.

In conclusion, the expression (\ref {Car2Bis}) can be rewritten as
\bee
v &=& F(u) + \sum_{\ell \in \N \, ,\ \ell \mbox { even}}  \frac {c^\ell}{\ell!} \frac {d^{(\ell-1)} }{du^{(\ell -1)}}\sum_{p=0}^\infty C^{(p)} \cos (pu) + \\ 
&& \ \ + \sum_{\ell \in \N \, ,\ \ell \mbox { odd}} \frac {c^\ell}{\ell !} \frac {d^{(\ell -1)} }{du^{(\ell -1)}}\sum_{p=0}^\infty D^{(p)} \sin (pu) \, , 
\eee
where the coefficients $ C^{(p)} $ and $ D^{(p)} $, which actually depend on the index $ \ell $ too, have been determined in the previous sections. \ Consequently the
series is reduced to a series of sine functions $\sin (pu)$ alone and its coefficients are generically denoted by $P$.

In section 20 the problem of determining the coefficient $ P $ of the generic term $ \sin (pu) $ is discussed; to this purpose he starts again from the series 
(\ref {Car2Bis}) distinguishing two contributions: the first contribution to the coefficient $P$ of $ \sin (pu) $ coming from the finite series 
$\sum_{\ell =1}^{p}  \frac {c^\ell}{\ell !} \frac {d^{(\ell -1)} \sin^\ell (u) F'(u)}{du^{(\ell -1)}}$ is denoted by $P'$; while the second contribution to the coefficient 
$P$ of $ \sin (pu) $ coming from the infinite series $\sum_{\ell =p+1}^{\infty} \frac {c^\ell }{\ell !} \frac {d^{(\ell -1)} \sin^\ell (u) 
F'(u)}{du^{(\ell -1)}}$ is denoted by $P''$. 

In order to find out these two values he makes use of the previously obtained results for the coefficients $ C^{(p)} $ and $ D^{(p)} $; in particular, in sections 21 
and 22 the following expression of $ P '$
\be
P' = \frac {2\alpha^p}{p} z = \frac {2c^p}{p(1+f)^p} z \, ,\label {Car6}
\ee
is given, where 
\bee 
z = e^{pf} \left ( 1 - \frac {1}{p!} \right ) \int_0^{pf} x^p e^{-x} dx \, . 
\eee

Sections from 23 to 30 are devoted to the calculation of the asymptotic behaviour of $ P '$ for large values ​​of $ p $ through a sequence of expansions in power series 
and additions by series. 

From section 31 he starts to compute the coefficient $ P '' $, where he makes some calculation errors. \ These errors alter the final result as pointed out by 
Jacobi. \ Initially, Carlini considers the case in which $ p = 2r $ is even; in the series $\sum_{\ell =2r+1}^{\infty}  \frac {c^\ell }{\ell !} \frac {d^{(\ell -1)} \sin^\ell (u) F'(u)}{du^{(\ell -1)}}$
he denotes by $ \Pi$ the coefficient of $ \sin (2ru) $ resulting from the $ 1^{\circ} $, $ 3^\circ $, $ 5^\circ $, etc. term of the series, and with $ \Pi'$ the 
coefficient of $ \sin (2ru) $ resulting from the $ 2^\circ $, $ 4^\circ $, $ 6^\circ $, etc. term of the series; consequently $ P '' = \Pi + \Pi'$. \ Since an expression of 
the term $D^{(\ell )}$ has been previously given then Carlini finally finds in section 31 that 
\be
\Pi = \frac {\left ( \frac 12 p \right )^p c^{p+1}}{p!} \sin^{p+1} (n) \left ( \alpha^p  z -  \alpha^{-p} z' \right ) \label {Car4_2}
\ee
where he sets\footnote {There are some sign errors; in fact the series are actually with alternate signs and the correct expressions are the following ones: 
\bee
z= \int \frac {\alpha^{-p}}{\sin^{p+2} n} dn - \frac {\left ( \frac 12 p \tan n \right )^2}{1\cdot (p+1)} \int \frac {\alpha^{-p}}{\sin^{p+4} n} dn +  
\frac {\left ( \frac 12 p \tan n \right )^4}{1\cdot 2\cdot (p+1)\cdot (p+2)} \int \frac {\alpha^{-p}}{\sin^{p+6} n} dn - ecc. \\ 
z'= \int \frac {\alpha^{p}}{\sin^{p+2} n} dn - \frac {\left ( \frac 12 p \tan n \right )^2}{1\cdot (p+1)} \int \frac {\alpha^{p}}{\sin^{p+4} n} dn + 
\frac {\left ( \frac 12 p \tan n \right )^4}{1\cdot 2\cdot (p+1)\cdot (p+2)} \int \frac {\alpha^{p}}{\sin^{p+6} n} dn - ecc.
\eee}
\bee
z&=& \int \frac {\alpha^{-p}}{\sin^{p+2} n} dn + \frac {\left ( \frac 12 p \tan n \right )^2}{1\cdot (p+1)} \int \frac {\alpha^{-p}}{\sin^{p+4} n} dn + \\ && \ \ + 
\frac {\left ( \frac 12 p \tan n \right )^4}{1\cdot 2\cdot (p+1)\cdot (p+2)} \int \frac {\alpha^{-p}}{\sin^{p+6} n} dn + ecc. \\ 
z'&=& \int \frac {\alpha^{p}}{\sin^{p+2} n} dn + \frac {\left ( \frac 12 p \tan n \right )^2}{1\cdot (p+1)} \int \frac {\alpha^{p}}{\sin^{p+4} n} dn + \\ && \ \ + 
\frac {\left ( \frac 12 p \tan n \right )^4}{1\cdot 2\cdot (p+1)\cdot (p+2)} \int \frac {\alpha^{p}}{\sin^{p+6} n} dn + ecc.
\eee

For the calculation of $ \Pi '$ he proceeds in a similar way making use of the results previously obtained regarding the coefficients $ C^{(p)} $. \ In particular, 
in section 32 he obtains that
\be
\Pi' &=& \frac {\left ( \frac 12 p \right )^{p+1} c^{p+1}}{(p+1)!} \sin^{p+2}(n) \left ( \alpha^p  z'' -  \alpha^{-p}  z''' \right ) \, , \label {Car4_3}
\ee
being\footnote {also in this case there are sign errors: the series have alternating signs as in the case of $ z $ and $ z '$.}
\bee
z'' &=& \int \frac {\alpha^{-p-1}}{\sin^{p+3} n} dn + \frac {\left ( \frac 12 p \tan n \right )^2} {1\cdot (p+2)} \int \frac {\alpha^{-p-1}}{\sin^{p+5} n} dn + \\ && \ \ + 
\frac {\left ( \frac 12 p \tan n \right )^4} {1\cdot 2\cdot (p+2)(p+3)} \int \frac {\alpha^{-p-1}}{\sin^{p+7} n} dn + ecc. \\ 
z''' &=& \int \frac {\alpha^{p+1}}{\sin^{p+3} n} dn + \frac {\left ( \frac 12 p \tan n \right )^2} {1\cdot (p+2)} \int \frac {\alpha^{p+1}}{\sin^{p+5} n} dn + \\ && \ \ + 
\frac {\left ( \frac 12 p \tan n \right )^4} {1\cdot 2\cdot (p+2)(p+3)} \int \frac {\alpha^{p+1}}{\sin^{p+7} n} dn + ecc. \\ 
\eee

If we set\footnote {Again, the correct series must has alternating signs; the first mistake is spreading along the paper!}
\be
s&=& s(x) := 1+ \frac {\left ( \frac 12 p x\right )^2}{1\cdot (p+1)\cdot \sin^{2} n}  + \frac {\left ( \frac 12 p x \right )^4}{1\cdot 2 \cdot (p+1)\cdot (p+2) 
\cdot \sin^{4} n} +  \nonumber \\ 
&& \ \ + \frac {\left ( \frac 12 p x \right )^6}{1\cdot 2 \cdot 3 \cdot (p+1)\cdot (p+2) \cdot (p+3) \cdot \sin^{6} n} + ecc. \label {Car4_4}
\ee
then it follows that\footnote {In fact, Carlini puts $x= \tan n$; however, in order to take into account the fact that the series has alternate signs the correct 
position should be $ x^2 = - \tan^2 (n) $, as Jacobi properly writes.}
\bee
z = \int \frac {\alpha^{-p}}{\sin^{p+2} n} s \left [ \tan (n) \right ] dn \, , 
\eee
and a similar expression can be given for the other three terms $z'$, $z''$ and $z'''$. \ Thus, the calculation of the series (\ref {Car4_4}) plays a crucial role 
in order to get the value of $\Pi$ and $\Pi'$. \ To this end Carlini applies a quite ingenious method; indeed, through some formal algebraic passages he finds 
that the sum of the series $ s = s (x) $ must solve the ordinary differential equation of the second order
\be
\frac {d^2s}{dx^2} + \frac {2p+1}{x} \frac {ds}{dx} = \frac {p^2}{\sin^2 (n)} s\, . 
\ee

In section 34 Carlini, by setting $ s = e^{\frac 12 p \int y dx} $, reduces this equation to the first-order equation 
\be
\frac {dy}{dx} + \frac 12 p y^2 + (2p+1) \frac {y}{x} = \frac {2p}{\sin^2 n} \, , 
\ee
recognising it as an equation already solved by Euler by means of definite integrals; in particular, recalling that $ \sin^2 n $ is equal to $ - \frac {1-c ^ 2}{c ^ 2} $, 
he finds that\footnote {Carlini is inaccurate in reporting the solution given by Euler. \ As Jacobi pointed out, the correct expression of the solution turns out to be
$s=C \int_0^{\frac {c}{\sqrt {1-c^2}}} \left ( \frac {c^2}{1-c^2} - t^2 \right )^{p-\frac 12 } \cos (p x t) dt $ where the numerical prefactor $C$ is given by  
$C= \frac {2\cdot 4 \cdots 2p}{1\cdot 3 \cdots (2p-1)} \frac {2}{\pi} \left ( \frac {1-c^2}{c^2} \right )^p$.
} 
\be
s= \int_0^{\frac {c}{\sqrt {1-c^2}}} \left ( \frac {c^2}{1-c^2} - t^2 \right )^{p-\frac 12 } \cos (p x^2 t) dt \, .\label {Car3}
\ee

According to Carlini, the formula (\ref {Car3}) has a purely theoretical meaning and it is not of great help. \ So Carlini tries to find an approximate value of the 
solution $ s $ when $ p $ takes very large values. \ Setting
\be
y = Y + \frac {Y'}{p} +  \frac {Y''}{p^2} +  ecc. 
\ee
Carlini finds a recursive procedure for determining the terms $ Y $, $ Y'$, $ Y'', \ldots $. \ More precisely, the dominant term $ Y $ in $ p $ is the solution to the 
equation
\bee
Y^2 +\frac {4Y}{x} = \frac {4}{\sin^2 n}\ \mbox { that is } \ Y= \frac 2x \left ( \sqrt {1+ \frac {x^2}{\sin^2 n}} - 1 \right )\, , 
\eee
while, in section 36, the subsequent terms are obtained through the recursive relations
\bee
\left ( Y + \frac 2x \right ) Y' &=& - \frac {dY}{dx} - \frac Yx \\ 
\left ( Y + \frac 2x \right ) Y'' &=& - \frac {dY'}{dx} - \frac {{Y'}^2}2 - \frac {Y'}{x} \\ 
\left ( Y + \frac 2x \right ) Y''' &=& - \frac {dY'}{dx} - Y' Y'' - 
\frac {Y''}{x} 
\eee
and so on, obtaining, in section 37, that
\bee
y = \frac 2x (g-1) - \frac 1p \frac {x}{g^2 \sin^n n} + \frac {1}{p^2} \left ( \frac {x}{g^3 \sin^2 n}- \frac 54 \frac {x^3}{g^5 \sin^4 n} \right ) + ecc. \, . 
\eee
where $ 1 + \frac {x^2} {\sin^2 n} = g^2 $. 

Since he had set $ s = e^{\frac 12 p \int y dx} $ then section 38 is dedicated to the calculation of the integral
\bee
\int y dx = 2(g-1) -2 \ln \left ( \frac {g+1}{2} \right ) - \frac 1p \ln g + \frac {1}{p^2} \left ( \frac 16 + \frac {1}{4g} - \frac {5}{12g^3} \right ) + ecc.  
\eee
where the integration constant is chosen so that the integral has value zero when $ x = 0 $, that is $ g = 1 $. 

Now, Carlini is ready to determine $ z $ and $ z '$ finding, in sections 39 and 40, the value, always understood as the dominant term in the limit of large $p$, 
for $ \Pi $. 

Similarly, in sections 41 and 42, the terms $ z '' $, $ z'''$ and the expression of the dominant term of $ \Pi' $ are determined. 

Having computed the values ​​of $ \Pi $ and $ \Pi' $, now he can find the expression of the dominant term of $ P '' = \Pi + \Pi '$ for large $p$; in 
particular, in section 43 it is obtained that
\be
P'' = \frac {4c^2}{A\sqrt {p\pi}} \frac {pf}{hg} \left ( \frac {c e^g}{1+g} \right )^p + ecc.\, ; \label {Eq10Bis}
\ee
where, being $ x = \tan n $, he sets \footnote {Recall that, as stated by Jacobi, the correct position is $x^2 =-\tan^2 n$.}
\bee
f &=& \sqrt {1-c^2} \, , \ g = \sqrt {1+c^2} \\
h &=& \left ( \frac {pc^2}{g+1} +2 \right )^2 + 2p \left ( \frac {pc^2}{g+1} +2 \right ) + p^2 c^2 \\ 
A &=& - \frac {p!e^p }{p^p\sqrt {p\pi }} \, . 
\eee
By virtue of this result and from formula (\ref {Car6}) he finally can obtain the value of the dominant term of $ P $, that is of the coefficient of $ \sin (pu) $, in 
the limit of large and even values of $ p $. 

In section 44 Carlini puts the question of the absolute convergence of the obtained series. \ Concerning the coefficients $P'$ he finds out that the associated series 
is a geometric series with common ratio
\bee
\frac {c}{1+f} e^f= \frac {c }{1+\sqrt {1-c^2}}e^{1-c^2}
\eee
which converges when this ratio is strictly less than unity, that is when 
\bee
\ln c <\ln \left ( 1+ \sqrt {1-c^2}\right ) - \sqrt {1-c^2}\, . 
\eee
Since this inequality holds true for any $ c <1 $, then it follows that for all the elliptical orbits the series, whose general term has coefficient determined by 
$ P '$, is convergent. 

On the other hand, the main term of the coefficients $P''$ is associated to a geometric series with common ratio
\be
\frac {c }{1+g}e^g = \frac {c }{1+\sqrt {1+c^2}} e^{ \sqrt {1+c^2}}\, ; \label {Err1} 
\ee
so with regard to $ P '' $ the series will be convergent provided that 
\bee 
\ln c < \ln \left ( 1+ \sqrt {1+c^2}\right ) - \sqrt {1+c^2}\, , 
\eee
which is satisfied when about $ c < 0.66$. 

In section 45 Carlini gives an example of calculation. 

In section 46 he quickly extends the results to the case in which $ p = (2r + 1) $ is odd. 

Finally, in section 47 he analyses the behaviour of the radius vector calculated as a series of cosine functions of multiples of 
the average anomaly finding that the $ p $-th coefficient of this expansion has expression given by
\be
Q = \pm \frac {2c^p e^{pg}}{2^p p A \sqrt {p\pi }} g \left ( \frac {2}{g+1} \right )^2 \, ; \label {RefQ}
\ee
consequently the convergence of the series will be guaranteed when the eccentricity $ c $ is less than the threshold value
such that 
\be
\ln c + \sqrt {1+c^2} = \ln 2\, ,\label {Arnold2}
\ee
that is when $c$ is less than the value $0.62$.

\section {Other mathematical works by Carlini} \label {AppB}

Carlini published more than 145 scientific articles and some of them were about mathematics. \ In this Appendix, following the cronological order, I briefly comment them showing that although they are less interesting than the Memoir by Carlini on Kepler's equation, they denote Carlini's remarkable predisposition for mathematics.

\subsection {Memoir by Carlini entitled {\it Sopra alcune funzioni esponenziali comprese nella formola $x^{x^n}$}} This Memoir was presented in Milan at the {\it Istituto del 
Regno Lombardo-Veneto} on the occasion of the meeting of January 7th, 1813 and then published in 1819 \cite {Esponenziale}. \ Therefore, it is one of Carlini's first works 
and, limited to his purely mathematical contributions, the first one. \ This 14-pages work, entirely in Italian, is written very clearly and it deserves to be 
analysed with some detail because addresses some interesting issues. \ Let us take up the whole introduction of the Memoir by Carlini:
\begin {itemize} 

\item [] {''After algebraic, exponential or logarithmic, and trigonometric functions, the simplest of its kind is the one represented by $ x^x $. \ The famous Giovanni Bernoulli, 
who was the first to consider it, discovered several different unique properties (see Johannis Bernoulli opera omnia, tom. I, pag. 184; tom. III, pag. 376), and he 
came to that very elegant series, which Leibniz liked so much,
\bee
\frac {1} {1 ^ 1} - \frac {1} {2 ^ 2} + \frac {1} {3 ^ 3} - \frac {1} {4 ^ 4} + etc. ,,
\eee
which is obtained by taking the integral of $ x^x dx $ from $ x = 0 $ to $ x = 1 $. \ Since I have spread the searches around this function somewhat further on, and around the more general $x^{x^n}$, I have come to find the sum of some equally 
divergent series, as convergent as the one referred to above. \ I have also examined the nature of the imaginary roots to the equation $ x^x = y $, which, as far as I know, they 
had not yet been studied by anyone. \ This is what I'm going to present in this brief Memoir.``}\footnote {\S 1, pg. 3 \cite {Esponenziale} ``Dopo le funzioni algebraiche, le esponenziali o lagaritmiche, e le 
circolari, la pi\`u semplice nel suo genere \`e quella rappresentata da $x^x$. \ Il celebre Giovanni Bernoulli, che primo prese a considerarla, ne scopr\`\i\ diverse 
diverse singolari propriet\`a (v. Johannis Bernoulli opera omnia, tom. I, pag. 184; tom. III, pag. 376), ed arriv\`o a quella serie elegantissima, e che tanto piacque 
al Leibnizio, 
\bee
\frac {1}{1^1} - \frac {1}{2^2} + \frac {1}{3^3} - \frac {1}{4^4} + ecc. \, ,
\eee
la quale si ottiene prendendo l'integrale di $x^x dx$ da $x=0$ ad $x=1$. \ Avendo io stese alquanto pi\`u avanti le ricerche intorno a questa funzione, ed intorno all'altra pi\`u generale $x^{x^n}$, sono giunto a trovare la somma di alcune serie 
altrettanto divergenti, quanto convergente \`e quella riferita di sopra. \ Ho pure esaminata la natura delle radici immaginarie dell'equazione $x^x =y$, le quali, per 
quanto ne so, non erano ancora state ricercate da alcuno. \ Ci\`o \`e quanto esporr\`o in questa breve Memoria.} 

\end {itemize}
In the twenty sections of the Memoir Carlini explores the properties of the solution $x=x(y)$ to the equation $ x^x = y $ (or its generalisation $ x^{x^n} = y $). \ Before 
analysing Carlini's work we recall that the equation $ x^x = y $ was already considered by Euler \cite {Eulero1} and, before him, by Lambert \cite {Lambert} in 1758. 

He recalls that the solution $ x $ to the equation $ x^x = y $ can be expressed through the power series
\be
x = 1+z+\sum_{n=2}^{+\infty} (-1)^{n-1} \frac {(n-1)^{(n-1)}}{n!} z^n  \, , \ \mbox { where } \ z = \ln y \, , \label {Eul1}
\ee
which turns out to be convergent when $ | z | <e^{- 1} $, by virtue of the (Stirling) formula used by the author \footnote {In fact, we know that this formula is 
correct when referred to the dominant term in the limit of large $n$.} $\frac {(n-1)^{(n-1)}}{n!} z^n = \frac {e^{n-1} z^n}{n \sqrt {2n\pi }}$. 

Then, Carlini pays attention to the question of the existence of the solution to the equation $ x \ln x = z $ when $ z <-e^{- 1} $. \ He properly acknowledged 
that the series (\ref {Eul1}) is divergent when $z<-e^{-1}$ and that the solutions cannot be real-valued, ​​but complex-valued of the form \footnote {We adopt the usual 
convention to denote the imaginary unit with $ i $, while Carlini adopts the notation $\sqrt {-1}$.}
\bee
x = -\frac {z \sin \alpha }{\alpha } e^{i\alpha }\, , 
\eee
where $ \alpha $ is a solution to the equation
\bee
\ln (-z) = \ln \left (\frac {\alpha}{\sin (\alpha )} \right ) - \frac {\alpha}{\sin (\alpha )} \, . 
\eee
It is interesting to note that Carlini recognises that the equation $ x \ln (x) = z $ can admit infinite solutions in the complex plane; in fact he states 
that

\begin {itemize}

\item [] {''...Since nothing restricts the value of $ z $ in these operations, it is clear that not only the imaginary roots to the equation $x\ln x =z$ will be 
found when, being $ z <- \frac 1e $, there is none real-valued, but still the imaginary roots are infinitely many, when the equation itself has one or
two real-valued solutions.``}\footnote {\S 5, pg. 5 \cite {Esponenziale} ``...Siccome niente limita in queste operazioni il valore di $z$, \`e chiaro che con esse si troveranno non solo le radici 
immaginarie dell'equazione $x\ln x =z$ quando, essendo $z<-\frac 1e$, non ve n'\`e alcuna di reale, ma ancora le radici immaginarie in numero infinite, quando l'equazione 
stessa ne ha una o due di reali.''}

\end {itemize}
It deserves to be emphasised the fact that the search for imaginary roots to equation $x^x =y$ is an original contribution given by Carlini; in fact, the first works, after that of Carlini, on 
the search for complex-valued roots to this kind of equations are those by L\'emeray \footnote {L\'emeray studied the real and, then, complex roots to the equation $ x = a^x $, related 
to the equation studied by Carlini, describing, in the case of complex solutions, the structure of the branches.} of the years 1896 \cite {L1} and 1897
\cite {L2}. \ However, Carlini's work did not have the diffusion that it deserved and his researches were ignored by L\'emeray.

Eventually, he begins the analysis of the integral $ \int x^{- \frac xr} dx $ with $ r> 0 $. \ By making the substitution $ z = -\frac xr \ln x$ and by integrating 
by parts Carlini finds the formal series
\bee
\int x^{-\frac xr} dx &=& \int e^z \frac {dx}{dz} dz = e^z \left ( \frac {dx}{dz} - \frac {d^2 x}{dz^2} + \frac {d^3 x}{dz^3} +  ecc. \right ) + C \\ 
&=& f(x) + C 
\eee
where
\bee
f(x) &=& x^{-\frac xr} \left [ - \frac rp + \frac {r^2}{xp^3} - \frac {r^3}{x^2} \left ( \frac {3}{p^5} + \frac {1}{p^4} \right ) + \frac {r^4}{x^3} 
\left ( \frac {15}{p^7} + \frac {10}{p^6} + \frac {2}{p^5} \right ) + \right. \\ 
&& \left. - \frac {r^5}{x^4} \left ( \frac {105}{p^9} + \frac {105}{p^8} + \frac {40}{p^7}+ \frac {6}{p^6} \right ) 
+ ecc. \right ] \, , 
\eee
$ C $ denotes a generic constant integration and $p= 1 + \ln x$. 

It should be noted that Carlini's reasoning has some flaws due to lack of mathematical accuracy. \ In fact, the integration by parts is correct when one
iterates such a procedure a finite number $ N $ of times and where the generic term $ etc. $ is meant to include the remainder integral 
$ \int (-1)^N e^z \frac {d^N x} {dz^N} dz $. \ If instead the integration by parts is intended to lead to infinite iterations and the term $ etc. $ only includes the 
terms $ (- 1)^ne^z \frac {d^nx} {dz^n} $ for $ n = 0, 1, 2, \cdots $ then we arrive to a paradox (which is a typical situation when one switches from a finite number of 
iterations to infinitely many without precautions!). \ It is not necessary to be too stern with Carlini: this Memoir was presented in 1813 when a rigorous approach 
to this kind of argument was still to come; indeed, Carlini's reference was still Euler's work where, with ingenious procedures, the calculation of  
divergent series was carried out. \ For example, Euler himself \cite {Eul} found that the divergent series $ \sum_{n = 0}^{\infty} (-1)^n n! $ has the integral 
$\int_0^1 \frac 1t e^{1-t^{-1}} dt$ as ``sum''.

In conclusion, in his Memoir Carlini addressed and solved three questions:
\begin {itemize}

\item [1.] The extension of the analysis concerning the problem $ x^x = y $ to the problem $ x^{x^n} = y $;
 
\item [2.] The calculation of the sum of the divergent series $1+ \sum_{n=1}^\infty (-n)^n$;
 
\item [3.] The search for complex solutions to the equation $x^x =y$.
 
\end {itemize}
In fact, the extension to the problem $ x^{x^n} = y $ does not seem really interesting, as well as the calculation of the divergent series 
$ 1 + \sum_{n = 1}^\infty (-n)^n $, also if compared to similar results \cite {Eul}. \ Conversely, the study of complex solutions appears to be innovative and 
worthy of greatest interest. \ The intuition that one may have infinite complex solutions anticipated an important result obtained by L\`emeray in 1897.

\subsection {Memoir by Carlini entitled {\it Legge dell'inserzione delle foglie nelle piante}} \label {Macchi} This is a brief memoir of popular nature published in 1837 in the 
magazine {\it Biblioteca Italiana} \cite {InserzioneFoglie}, at that time directed by Carlini himself. \ In this work Carlini summarises the recent researches of some foreign authors 
(for example the note by Alexander Braun of the year 1831 \cite {Braun} or the note by Karl Friedrich Schimper of the year 1836 \cite {Schimper}) dedicated to phyllotaxis, which is 
that branch of botany devoted to the study and determination of the order by which the various botanical entities (leaves, flowers, etc.) are distributed in space, 
giving plants a geometric structure. \ In particular, the results presented in the year 1837 by Louis and Auguste Bravais in a memoir at the {\it Acad\'emie royale des 
Sciences} in Paris \cite {Bravais} are detailed summarised. \ In this Memoir Carlini explains that the leaves, instead of being separated from each other by an angle equal to 
$ \frac 13 $, or to $ \frac 25 $, or to $ \frac  38$, etc., of the circumference, are separated by a constant angle, incommensurable with the circumference itself, and, 
consequently, such that two leaves can never be located exactly along the same vertical direction. \ In fact, to this end the divergent numerical sequence $1$, $2$, $3$, 
$5$, $8$, etc., is introduced, where 
\begin {itemize} 

\item [] {''...each of which equals the sum of the previous two...``}\footnote {Pg. 287 \cite {InserzioneFoglie} ``...ciascuno de' quali eguaglia la somma dei due precedenti...''}

\end {itemize}
and it has been concluded that the incommensurable angle that one wants to determine is equivalent to $ 2 \pi \lambda $, where $ \lambda = \frac {-1+ \sqrt {5}} {2} $ is the 
positive root of the trinomial $ x ^ 2 + x-1 $. \ In short, mathematical analysis concerns those that were later called \emph {Fibonacci sequence} and 
\emph {golden ratio}. \footnote {The Fibonacci sequence was introduced by Fibonacci in his \emph {Liber Abaci}, published as a manuscript in 1202. \ However, it must 
be said that the first printed version, which had a certain diffusion, appeared only in 1857 by the mathematician and science historian Baldassarre Boncompagni; in 1877 
the French mathematician \'Edouard Lucas shown profound results on Fibonacci numbers and he indicated that such a sequence is found for the first time in the Liber Abaci 
and begun to call it Fibonacci sequence, a name that was later universally adopted. \ Concerning the origin of the name \emph {golden ratio} we refer to \cite {Golden1,Golden2}.}

\subsection {Memoir by Carlini entitled {\it Descrizione d'una macchinetta che serve a risolvere il problema di Keplero, ossia a trovare l'anomalia eccentrica data 
l'anomalia media, qualunque sia l'eccentricit\`a}} This is a brief memoir of the report presented by Carlini at the meeting of the {\it Istituto Lombardo Veneto} of August 5, 1852, 
published then in 1853 \cite {macchinetta}. \ In this report Carlini returns to treating Kepler's equation by proposing a numerical procedure using a kind of an analogue computer for its solution even 
in the presence of not too small eccentricity values. \ As an alternative to a graphic method for finding the solution, Carlini designs a ``macchinetta'' (that is a small machine)  
that allows us to implement a solution method of Kepler's equation by successive approximations. \ With modern language Carlini deals with Kepler's equation $ x + e \sin x = z $, 
where $ 0 \le e \le 1 $, as a fixed point problem $ x = f (x) $ where $ f (x) = z-e \sin x $; using the machine he finds an approximate $ x_0 $ value of the solution, after 
which the iterative process is activated $ x_{i + 1} = f (x_i) $ starting from the index $ i = 0 $ and arriving in a few steps to an acceptable approximation of the solution. \ The  
description of the machine operation is given in the memoir and a sketch of the device is also drawn, but both the description of the procedure and the device 
diagram are not very clear, and in any case the device does not seem to have had a diffusion of any importance.

\subsection {Memoir by Carlini entitled {\it Sulle propriet\`a delle funzioni algebriche conjugate}} This Memoir, written in Italian language, was presented to the 
{\it Sitzungsberichte der Akademie der Wissenschaften} (Wien) in July 1854 \cite {FunzioniConiugate}, so at a mature age. \ In this Memoir, of purely algebraic-analytic 
kind, the author explores the properties of the two conjugated algebraic functions
\bee
f(x) = x - \frac 1x \ \mbox { and } \ \varphi (x) = x + \frac 1x \, . 
\eee
After reviewing the main properties, Carlini considered the application to the numerical calculation of real solutions to second and third degree equations. \ For 
example, the general equation $ p ^ 2-ap = b $, $ b> 0 $, can be placed in the form $ f (x) = \frac {a}{\sqrt {b}} $ through the replacement $ p = x \sqrt {b} $. \ Thus, 
by making use of the numerical tables of the functions $ f (x) $ and $ \varphi (x) $, placed in the appendix of the Memoir, the numerical solution to the original equation easily 
follows. \ More interesting is the application to third degree equations of the type $ p ^ 3 + ap + b $, $ a> 0 $; if one sets $ p = q \sqrt {a / 3} $
then the equation takes the form $ q^3 + 3p = (3 / a)^{3/2} b $. \ On the other hand, with the property $ f (x^3) = [f (x)]^3 + 3f (x) $, the latter equation can be
rewritten as $ f (y) = (3 / a)^{3/2} b $, where $ y = q^3 $, which admits immediate numerical solution using the aforementioned tables.


\begin{thebibliography}{99}

\bibitem {A} V.I.Arnold, {\it Huygens and Barrow, Newton and Hooke. \ Pioneers in mathematical analysis and catastrophe theory from evolvents to quasicrystals}, 
Birkh\"auser Basel (1990).

\bibitem {Golden1} H.Becker, {\it An even earlier (1717) usage of the expression ``golden section''}, Historia Mathematica {\bf 49} 82-83 (2019).

\bibitem {BM} M.V.Berry, and K.E.Mount, {\it Semiclassical approximation in wave mechanics}, Rep. Prog. Phys. {\bf 35} 315-389 (1972).

\bibitem {Bessel} F.W.Bessel, {\it Analytische Aufl\"osung der Keplerschen Aufgabe}, Abh. Preuss. Akad. Wiss. Berlin, {\bf XXV} 49-55 (1819).

\bibitem {Bottazzini} U.Bottazzini, and J.Gray, {\it Hidden Harmony-Geometric Fantasies. \ The rise of complex function theory}, Springer Verlag New York (2013).

\bibitem {Braun} A.Braun, {\it Vergleichende Untersuchung \"uber die Ordnung der Schuppen an den Tannenzapfen als Einleitung zur
Untersuchung der Blattstellung \"uberhaupt}, Nova Acta Ph. Medical Academy of Cesar Leop Carolina Nat Curiosorum {\bf 15}, 195-402 (1831).

\bibitem {Bravais} L.Bravais, and A.Bravais, {\it Essai sur la disposition des feuilles curvis\'eri\'ees}, Annales des Sciences Naturelles Botanique, {\bf 7} 42-110 (1837).

\bibitem {Carlini1} F.Carlini, {\it Ricerche sulla convergenza della serie che serve alla soluzione del problema di Keplero}, Dall'Imp. Regia Stamperia Milano 1-48 (1817).

\bibitem {Esponenziale} F.Carlini, {\it Sopra alcune funzioni esponenziali comprese nella formola $x^{x^n}$}, Memorie dell'I.R. Istituto del Regno Lombardo-Veneto 167-178 (1819).

\bibitem {InserzioneFoglie} F.Carlini, {\it Legge dell'inserzione delle foglie nelle piante}, Biblioteca Italiana (Giornale di letteratura scienze ed arti) {\bf 86} 296-288 (1837).

\bibitem {macchinetta} F.Carlini, {\it Descrizione d'una macchinetta che serve a risolvere il problema di Keplero, ossia a trovare l'anomalia eccentrica data 
l'anomalia media, qualunque sia l'eccentricit\`a (con una tavola)}, Giornale dell'I.R. Istituto Lombardo di Scienze, Lettere ed Arti {\bf 5} 438-445 (1853).

\bibitem {FunzioniConiugate} F.Carlini, {\it Sulle propriet\`a delle funzioni algebriche conjugate. \ Memoria presentata all'I.R. Academia delle Scienze di 
Vienna dal Membro effettivo Francesco Carlini (con una tavola)}, Sitzungsberichte der Kaiserlichen Akademie der Wissenschaften mathematisch-naturwissenschaftliche 
Classe, {\bf 13} 357-375 (1854).

\bibitem {Colwell} P.Colwell, {\it Solving Kepler's equation over three centuries}, Willmann-Bell Inc. (1993).

\bibitem {Dutka} J.Dutka, {\it On the early history of Bessel functions}, Archive for History of Exact Sciences {\bf 49} 105-134 (1995).

\bibitem {Eul} L.Euler, {\it De seriebus divergentibus}, Novi Commentarii Academiae Scientiarum Petropolitanae {\bf 5} 205-237 (1760).

\bibitem {Eulero1} L.Euler, {\it De serie lambertina plurimusque eius insignibus propriatetibus}, Acta Academiae Scientiarum Imperialis Petropolitanae {\bf 2} 29–51 (1779).

\bibitem {Finch} S.R.Finch, {\it Mathematical constants}, Cambridge University Press (2003).

\bibitem {Froman} N.Fr\"oman, and P.O.Fr\"oman, {\it Physical Problems Solved by the Phase-Integral Method}, Cambridge University Press (2002).

\bibitem {FF2} N.Fr\"oman, and P.O.Fr\"oman, {\it Stark effect in a hydrogenic atom or ion: treated by the phase-integral method}, Imperial College (2008).

\bibitem {Giacobbe} G.C.Giacobbe, {\it Francesco Carlini}, Dizionario Biografico degli Italiani, Ed. Treccani {\bf 20} (1977).

\bibitem {Golden2} R.Herz-Fischler, {\it An early usage of the expression ``golden section''}, Historia Mathematica {\bf 49} 80-81 (2019).

\bibitem {Jacobi1} C.G.J.Jacobi, {\it \"Uber die ann\"ahernde Bestimmung sehr entfernter Glieder in der Entwickelung der elliptischen Coordinaten nebst 
einer Ausdehnung der Laplaceschen Methode zur Bestimmung der Functionen gerader Zahlen},  Schumacher Astronomische Nachrichten {\bf 665} 257-270 (1849).

\bibitem {Jacobi2} C.G.J.Jacobi, {\it Untersuchungen \"uber die Convergenz der Reihe, durch welche das Kepleresche Problem gel\"ost 
wird. \ Von Franz Carlini}, Schumacher Astronomische Nachrichten {\bf 709-712} 197-254 (1850).

\bibitem {Lac1} S.F.Lacroix, {\it Trait\'e du calcul diff\'erentiel et du calcul int\'egral. Tome second}, Paris (1798).

\bibitem {Lac2} S.F.Lacroix, {\it Trait\'e du calcul diff\'erentiel et du calcul int\'egral. Seconde \'edition, revue et augment\'ee. Tome second}, Paris (1814).

\bibitem {Lagrange} J.L.Lagrange, {\it Sur le probl\'eme de Kepler}, Mem. de l'Acad. des Sci. Berlin, {\bf XXV} 204-233 (1771).

\bibitem {Lag} J.L.Lagrange, {\it M\'ecanique analytique}, Paris (1788).

\bibitem {Lambert} J.H.Lambert, {\it Observationes variae in mathesin puram}, Acta Helvetica, physico-mathematico-anatomico-botanico-medica, {\bf 3} 128-168 Basen (1758).

\bibitem {Laplace0} P.S.Laplace, {\it Sur le d\'eveloppement de l'anomalie vraie et du rayon vecteur elliptique, en s\'eries ordonn\'ees suivant 
les puissances de l'excentricit\'e}, M\'emories del l'Acad\'emie Royale des Sciences Paris {\bf 6}, 61-80 (1823).

\bibitem {Laplace} P.S.Laplace, {\it Suppl\'ement au $5^\circ$ volume du trait\'e de M\'ecanique C\'eleste}, Paris (1827).

\bibitem {L1} E.M.L\'emeray, {\it Sur les racines de l'\'equation $x=a^x$}, Nouvelles Annales de Math\'ematiques III series, tome {\bf 15} 548-556 (1896).

\bibitem {L2} E.M.L\'emeray, {\it Sur les racines de l'\'equation $x=a^x$. Racines imaginaires}, Nouvelles Annales de Math\'ematiques III series, tome {\bf 16} 54-61 (1897).

\bibitem {Moulton} F.R.Moulton, {\it An introduction to Celestial Mechanics}, The MacMillan Company, New York (1914).

\bibitem {Scheibner} W.Scheibner, {\it On the asymptotic values of the coefficients in the development of any power of the radius-vector 
according to the mean anomaly}, Astronomical Journal {\bf 95}, 177-182 (1856).

\bibitem {Schiaparelli1} G.V.Schiaparelli, {\it Notizie sulla vita e sugli studj di Francesco Carlini}, Atti del R. 
Istituto Lombardo di Scienze e Lettere, {\bf III} 281-292 (1862).

\bibitem {Schiaparelli2} G.V.Schiaparelli, {\it Notizie sull'osservatorio di Brera in Milano}, in {\it Mediolanum}, 
Casa Editrice Dottor Francesco Vallardi, Bologna, Milano, Napoli (1881).

\bibitem {Schimper} K.F.Schimper C.F., {\it Geometrische Anordnung der um eine Axe periferischen Blattgebilde}, Verhandl Schweiz Naturf Ges, {\bf 21} 113-117 (1836).

\bibitem {Schlissel} A.Schlissel, {\it The development of Asymptotic Solutions of Linear Ordinary Differential Equations, 1817-1920}, Archive for History of 
Exact Sciences {\bf 16} 307-378 (1977).

\bibitem {Taff} L.G.Taff, {\it Celestial Mechanics: a computational guide for the practitioner}, Wiley-Interscience (New York) (1985).

\bibitem {TT} G.Tagliaferri, and P.Tucci, {\it Carlini and Plana on the theory of the moon and their dispute with Laplace}, Annals of Science {\bf 56} 221-269 (1999).

\bibitem {Tournes} D.Tourn\`es, {\it L'int\'egration approch\'ee des \'equations diff\'erentielles ordinaires (1671-1914)}, Thesis (1996).

\bibitem {Tre} J.Trembley, {\it R\'eflexions sur l'usage des m\'ethodes d'approximation dans l'int\'egration des \'equations diff\'erentielles}, M\'emoires de 
l'Acad\'emie Royale des Sciences et Belles-Lettres Berlin 387-422 (1792).  

\bibitem {W} G.N.Watson, {\it A treatise on the Theory of Bessel functions}, Cambridge University Press (second edition) (1966).

\bibitem {WW} E.T.Whittaker, and G.N.Watson, {\it A course of modern analysis}, Cambridge University Press (third edition) (1920).

\end{thebibliography}
\end{document}